\documentclass[11pt]{article}

\usepackage{amssymb,amsmath,amsfonts}
\usepackage{enumerate}
\usepackage[margin=1in]{geometry}
\usepackage{algorithm}
\usepackage{graphicx}
\usepackage{algpseudocode}
\usepackage{booktabs}
\usepackage{multirow}
\usepackage{tabulary}
\usepackage{longtable}
\usepackage{wrapfig}
\usepackage{mathtools}
\usepackage{subfig}
\usepackage{hyperref}
\usepackage{url}
\usepackage{breakurl}
\usepackage[normalem]{ulem}

\newtheorem{thm}{Theorem}[section]

\newtheorem{rem}[thm]{Remark}

\newcommand\Rmn[2]{\R^{#1 \times #2}}

\def\R{\mathbb R}
\def\C{\mathbb C}
\def\bigO{\mathcal O}

\def\thalf{\tfrac{1}{2}}
\def\transsym{\mathsf{T}}
\newcommand{\tp}[1][]{^{{#1}\transsym}}

\def\initgain{\gamma}
\def\zsa{\alpha_{D}}
\def\zsamargin{\alpha_\mathrm{c}}

\makeatletter
\let\c@equation\c@thm
\makeatother
\numberwithin{equation}{section}

\title{A Comparison of Nonsmooth, Nonconvex, Constrained Optimization Solvers for the Design of Time-Delay Compensators}

\author{Vyacheslav Kungurtsev\thanks{
Department of Computer Science, Faculty of Electrical Engineering, Czech Technical University in Prague, 
\href{mailto:vyacheslav.kungurtsev@fel.cvut.cz}{\texttt{vyacheslav.kungurtsev@fel.cvut.cz}}
}
\and
Tim Mitchell\thanks{
Max Planck Institute for Dynamics of Complex Technical Systems Magdeburg, 
\href{mailto:mitchell@mpi-magdeburg.mpg.de}{\texttt{mitchell@mpi-magdeburg.mpg.de}}
}
\and
Tom\'a{\v s} Vyhl\'idal\thanks{
Department of Instrumentation and Control Engineering, Faculty of Mechanical Engineering, and Czech Institute of Informatics, Robotics and Cybernetics, Czech Technical University in Prague,
\href{mailto:tomas.vyhlidal@fs.cvut.cz}{\texttt{tomas.vyhlidal@fs.cvut.cz}}
}
}

\date{December 30, 2018}

\begin{document}
\maketitle

\begin{abstract}
We present a detailed set of performance comparisons of two state-of-the-art solvers for 
the application of designing time-delay compensators, an important problem in the field of robust control.
Formulating such robust control mechanics as constrained optimization problems often involves objective and constraint functions 
that are both nonconvex and nonsmooth, both of which present significant challenges to many solvers and their end-users 
hoping to obtain good solutions to these problems.
In our particular engineering task, the main difficulty in the optimization arises in a nonsmooth and nonconvex stability constraint, 
which states that the infinite spectrum of zeros of the so-called shaper should remain in the open left half-plane.
To perform our evaluation, we make use $\beta$-relative minimization profiles, recently introduced visualization tools
that are particularly suited for benchmarking solvers on nonsmooth, nonconvex, constrained optimization problems.
Furthermore, we also introduce new visualization tools, called Global-Local Profiles, which for a given problem and a fixed computational budget,
assess the tradeoffs of distributing the budget over few or many starting points, with the former getting more budget per point and latter less.
\end{abstract}

\section{Introduction}
\label{sec:intro}
Consider the following general constrained optimization problem
\begin{equation}
	\label{eq:nlp}
	\min_{x\in\R^n} \, f(x) \quad \text{s.t.} \quad 
	\begin{cases}
	c_j(x)\le 0 \qquad j\in\{1,\ldots,m\} \\
	c_j(x) = 0 \qquad j\in\{m+1,\ldots,q\}
	\end{cases},
\end{equation}
where the objective function $f : \R^n \rightarrow \R$ and each individual (in)equality constraint
function $c_j : \R^n \rightarrow \R$ are all continuous but may possibly be nonsmooth.
We assume that the functions are still continuously differentiable almost everywhere;
as such, the nonsmoothness of $f$ and each $c_j$ is limited to a set of measure zero in $\R^n$.  
To date, much of the literature for \emph{nonsmooth} optimization has focused on unconstrained minimization,
without any $c_j$ constraint functions; e.g. see \cite{BagKM14}.
However, a couple of new solvers specifically intended for nonsmooth, constrained optimization have recently 
appeared, opening up new possibilities for reliably solving actual applications that are modeled by \eqref{eq:nlp} 
in practice, despite the challenging properties of the objective and/or constraint functions.

In 2012, \cite{CurO12} proposed combining gradient sampling (GS), 
a nonsmooth optimization technique for unconstrained minimization \cite{BurLO05},
with sequential quadratic programming (SQP) to additionally handle (possibly nonsmooth) constraints.  
The resulting algorithm, simply called SQP-GS, was shown to have provable convergence to stationary points of \eqref{eq:nlp}, 
assuming that all the functions $f$ and $c_j$ are locally Lipschitz and that such stationary points are sufficiently \emph{calm} 
(for details on calmness, see \cite{Bur91,Bur91a}).  
On the other hand, the gradient-sampling procedure employed within SQP-GS requires evaluating $\bigO(n)$ gradients 
(for each of the particular $f$ and $c_j$ functions that are nonsmooth) at every iteration,
which can make SQP-GS a computationally expensive method, especially if the cost to compute these functions is nontrivial.  
This potential high cost is further exacerbated by the fact that nonsmooth optimization methods generally have 
slow convergence rates (i.e. typically linear),
due to the  lack of exploitable smoothly varying curvature at nonsmooth minimizers,
and can thus incur many iterations before satisfactorily meeting stationarity conditions.

More recently, in 2017, a philosophically different algorithm \cite{CurMO17} for solving \eqref{eq:nlp} was proposed,
i.e. one that does not rely on gradient sampling and 
which instead typically only requires $\bigO(1)$ gradient evaluations on average per iteration,
thus potentially making it significantly faster than SQP-GS.
Inspired by the observed efficiency and reliability of quasi-Newton methods for nonsmooth unconstrained optimization \cite{LewO13}, 
specifically Broyden-Fletcher-Goldfarb-Shannon (BFGS) updating, 
the authors of \cite{CurMO17} proposed a new BFGS-SQP method as a fast and reliable alternative to SQP-GS. 
This BFGS-SQP algorithm is implemented in a code called
GRANSO: GRadient-based Algorithm for Non-Smooth Optimization \cite{granso}.
Although there are currently no theoretical convergence guarantees for the underlying BFGS-SQP method,
it was shown that GRANSO can be a competitive solver in practice compared to SQP-GS \cite[Section 6]{CurMO17}.
On a large test set of challenging nonsmooth, nonconvex, constrained optimization problems \cite{psarnot},
half of which were not even locally Lipschitz,
GRANSO often fared well relative to SQP-GS in terms of minimization quality and finding the feasible region, 
all while being many times faster than SQP-GS
(typically over an order of magnitude faster, which, for those problems, was in line with the expectation of being faster by a factor of $\bigO(n)$ 
as the dominant cost in the running times was the function/gradient evaluations).
Additionally, the ``off-the-shelf" solvers (i.e., those not specifically developed for nonsmooth problems)
included as baselines in the evaluation done in \cite{CurMO17} did not perform well at all,
highlighting the necessity for purpose-built solvers for nonsmooth optimization.

Returning to our specific motivation of using nonsmooth, constrained optimization solvers to design time-delay compensators, 
so-called \emph{input shapers}, the overall goal is to shape the movements 
of a main-body system (e.g. a crane trolley) in order to damp
the oscillatory modes of a connected flexible subsystem (e.g. a payload suspended from a crane trolley).
Cast as an optimization problem, the design of such input shapers in this paper takes the following general form:
\begin{equation}
\label{eq:qp_nsineq}
\min_{x\in\R^n} \, x\tp Q x \quad \text{s.t.} \quad 
	\begin{cases} 
	A_1 x & \le b_1 \\
	A_2 x & = b_2 \\
	s(x) 	& \le 0
	\end{cases},
\end{equation}
where $Q \in \Rmn{n}{n}$ is positive definite, $A_1 \in \Rmn{m}{n}$, $A_2 \in \Rmn{p}{n}$, 
and $s : \R^n \to \R$ is a continuous but nonsmooth and nonconvex function that is possibly also not locally Lipschitz.
The quadratic objective function models desired performance characteristics in the frequency domain,
while the linear constraints specify necessary structural requirements.
The nonlinear and nonsmooth inequality constraint $s$ imposes a stability constraint on the shaper, 
namely that its spectrum, 
which is infinite for the given class of time-delay systems, must entirely reside in the open left half-plane.
Note however that \eqref{eq:qp_nsineq} is applicable to a wide range of
related problems in control system theory 
where fast and robustly stable transient characteristics must be balanced.

Without explicitly enforcing stability of the shaper, the 
resulting shaper may have unstable modes which could destabilize the closed-loop system \cite{VyhHKetal16}.
Consequently, in \cite{KunPVetal17}, designing a \emph{stable} input shaper via optimizing a specific form of \eqref{eq:qp_nsineq} 
with constraint $s$ enforcing stability was proposed, 
which the authors noted was the first use of optimization to stabilize an input shaper via moving its zeros into the open left half-plane.
However, attempting to solve \eqref{eq:qp_nsineq} with SQP-GS, the only solver specifically designed for 
nonsmooth, nonconvex, constrained optimization that was available to the authors at the time of those initial experiments, was not 
a straightforward matter.  First, since the constraint function $s$ may not be locally Lipschitz, 
the convergence guarantees for SQP-GS might not even apply for \eqref{eq:qp_nsineq}.
Second, these convergence guarantees, even if they did apply, do not say anything
about whether or not a resulting computed solution will be feasible;
in other words, the convergence guarantees do not preclude SQP-GS from converging to an infeasible stationary point.
Third, as the authors of \cite[Section 4, p. 13327] {KunPVetal17} noted, 
SQP-GS was indeed quite slow in practice.
While a stable input shaper was successfully obtained by optimizing a specific example of \eqref{eq:qp_nsineq} via SQP-GS 
\cite[Section 4]{KunPVetal17}, 
the difficulties of doing so were beyond the scope of the article and thus were not reported in detail, e.g.
how many starting points were required in order to find at least one stable input shaper and 
exactly how long did it take SQP-GS to compute these shapers.
These are relevant issues for engineers actually attempting to employ this new methodology
and with the recent availability of GRANSO, which might potentially solve \eqref{eq:qp_nsineq}
significantly faster than SQP-GS, there is now also the question of which solver to use in this or similar applications.
As such, the questions of what is the current best practice for solving input shaper design problems of the form \eqref{eq:qp_nsineq}, 
and at what expense, remain.   In this paper, we aim to address these concerns.
Furthermore, 
as alluded to earlier,
our evaluation of SQP-GS and GRANSO can potentially be informative 
for other related robust control design tasks as these often involve nonsmoothness
arising from stability constraints and/or max functions.

Of course, comparing the performance of competing solvers on a 
nonsmooth, nonconvex, constrained optimization problem is a multi-faceted affair.  
First, just the ability of a method to consistently find feasible solutions is paramount, 
particularly for nonlinear constraints since feasible points may not be knowable a priori and 
methods may not always find the feasible region for every starting point.
Given feasibility, of second importance is the quality of objective function minimization, that is, 
how far has the objective function been reduced.  
With nonconvexity, this is not a well-defined goal since methods could converge to different local minimizers.
Furthermore, with a possibly non-locally-Lipschitz function present in our optimization problem, 
there are no guarantees that the methods will not stagnate before reaching stationary points.
However, from a user's perspective, whichever feasible answer minimizes the objective function the most 
would be preferred.
Finally, the computational effort necessary to meet feasibility requirements and sufficiently minimize the objective
must also be considered.  Beyond different per-iteration costs that methods may have, this metric of interest is also
complicated by the facts that (a) the rate of progress of any given method may be highly variable during optimization 
and (b) that reasonable computational budgets of users cannot typically be reliably estimated or predicted by others.
As such, to increase the utility of benchmarks for users, it can be imperative to consider multiple budgets.
To handle such complicated benchmarking objectives, where the three goals of 
finding feasibility, objective minimization, and computational efficiency 
are also often at odds with one another, $\beta$-Relative Minimization Profiles ($\beta$-RMPs or just RMPs)
were proposed as new visualization tools to concisely address all these facets \cite[Section 5]{CurMO17}.
To help assess what is the best practice for solving our input shaper design problem, we employ
$\beta$-RMP benchmarking to compare SQP-GS and GRANSO.

Finally, we propose new visualization tools, called \emph{Global-Local Profiles (GL-Profiles)}, 
to help resolve a particular quandary when it comes to nonconvex and/or nonsmooth optimization, 
where different local minima may be found and/or where convergence might be quite slow. 
Given a fixed computational budget, the question is then how to best utilize that budget in 
order to obtain good solutions to such a problem.
One could devote the budget to runs from just one or a few starting points, 
to allow local minimizers to be obtained precisely (i.e. reaching stationarity conditions).
Alternatively, one could distribute the budget over many starting points, where the much smaller budget per starting point means
that minimizers will typically not be found so accurately but nevertheless, the best minimizer 
amongst those obtained 
could be much better than what would be found by only using few starting points.
GL-Profiles present this tradeoff in a concise and intuitive manner and describe not only the relative performance 
differences of the solvers but also
speak to the ease or difficulty of finding good solutions that is inherent to a particular nonconvex/nonsmooth problem of interest itself.
  
This paper is organized as follows.  In Section~\ref{sec:shaper},
we provide background and details on the input shaper design optimization problem with the nonsmooth stability constraint.
In Section~\ref{sec:solvers}, we first give a brief overview of the SQP-GS algorithm and then the BFGS-SQP method employed within GRANSO.
We describe our experimental setup for designing stable input shapers in Section~\ref{sec:experiments}, 
while also giving a preliminary analysis of our results.
In Section~\ref{sec:rmps}, we provide a primer on $\beta$-RMPs for the convenience of the reader and 
then do a $\beta$-RMP benchmark comparing SQP-GS and GRANSO for our input shaper design problem.
To address the question of how to best utilize a fixed computational budget for nonconvex/nonsmooth optimization problems,
we first introduce our new GL-Profiles and then do a GL-Profile benchmark 
to compare both methods with this new complementary perspective in Section~\ref{sec:glprofiles}.
Final remarks are given in Section~\ref{sec:conclusion}.

\section{The input shaper design problem using distributed delays}
\label{sec:shaper}
Input shaping is a well-known technique for damping undesirable oscillatory modes of mechanical systems,
to intelligently control the motions of a main-body such that when it comes to a stop, 
the flexible subsystem also returns to a stand still as quickly as possible.
Input shaping via the use of time-delay
filters was first proposed and implemented by Smith~\cite{smith1957posicast} in the 1950s 
and then further extended to robust and multi-mode design 
tasks~\cite{singer1990preshaping,singhose1994residual,tuttle1994zero} in the 1990s.
In addition to the aforementioned typical crane problem \cite{vaughan2010control,singhose2008input},
input shapers have also been used in many other applications, such as position control of flexible manipulators \cite{pereira2009adaptive} and 
vibration suppression in robotics \cite{park2006design}. 
For an extensive review on input shaping over the last 50 years, see \cite{singhose2009command}.

The specific control application represented by \eqref{eq:qp_nsineq} that we consider here 
arises from the use of a \emph{distributed-delay shaper} \cite{vyhlidal2013signal,vyhlidal2015parameterization}, 
where the inverse of the designed shaper is incorporated into the feedback loop of a control system,
a technique proposed in \cite{vyhlidal2015feedback} (see also \cite{vyhlidal2016inverse}).
The benefit of such an approach is that it 
extends the applicability of input shaping to compensating the flexible modes induced by both reference signal and disturbances.
In~\cite{vyhlidal2015parameterization}, an optimization-based approach for designing such an input shaper
was proposed with the form of~\eqref{eq:qp_nsineq} but without the constraint $s(x)<0$ and
where, for robustness, the objective function was a least-squares measure of the residual vibrations across
a small frequency range around the desired oscillation to be compensated.
However, by including the inverse of the shaper in the feedback loop, 
zeros of the shaper are likely to be ``turned into" poles of the closed-loop system,
as addressed in \cite{vyhlidal2015feedback}. 
As closed-loop stability is required, it thus makes sense to add the condition
that all zeros of the input shaper must reside in the open left half-plane, which 
will be modeled by the $s(x) < 0$ constraint in \eqref{eq:qp_nsineq}.
As mentioned in the introduction, a preliminary optimization-based approach using SQP-GS 
to enforce this nonsmooth stability constraint was considered in \cite{KunPVetal17}.

In more detail, we first describe an input shaper with a distributed delay.  
Given $T > 0$, a constant specifying the longest delay, let $\{\tau_1,\ldots,\tau_n\}$ be the set of
equally distributed delays in $[0,T]$, indexed in increasing order, with $\tau_1 = 0$ and $\tau_n = T$.
The specific input shaper with delays proposed in \cite{vyhlidal2015parameterization} is characterized by the equation
\begin{equation}
\label{eq:shaper_dyn}
y(t)=\gamma u(t)+\int^{T}_{0}{u(t-\mu)dh(\mu)},
\end{equation}
where $u$ and $y$ are respectively the shaper input and output, $\gamma \in (0,1)$ is the gain parameter, and the distribution of the delays is prescribed by the non-decreasing function $h(\mu)$, with length $T$
(though again, we here will only focus on equally-distributed delays).
Considering the step-wise distribution of the delays, the
corresponding transfer function $D: \C \to \C$ of \eqref{eq:shaper_dyn} has the following form:
\begin{equation}
\label{shpDZVTF}
D(s)=\initgain+\frac{1}{s} \sum_{k=1}^{n}{x}_k{\rm e}^{-s\tau_k},
\end{equation}
where $x_k$ is the respective gain for each delay spike with delay $\tau_k$.
The set of $x_k$  values comprise our set of optimization variables.

In the classical feed-forward interconnection of the shaper $D(s)$ and the flexible subsystem $F(s)$, which acts as $D(s)F(s)$, 
the mode compensation takes place via cancellation of the oscillatory pole of $F(s)$ by a zero of the shaper $D(s)$. 
The purpose of including the inverse shaper is to pre-compensate the oscillatory modes of the flexible subsystem 
$F(s)$ when controlling the main-body system $G(s)$ via the controller $C(s)$.  In contrast to the classical feed-forward 
application of the shaper $D(s)F(s)$, 
incorporating the inverse of the shaper allows flexible modes that are excited by disturbances $d$ to also be compensated. This feature results from the corresponding transfer functions from the system reference $w$ to the system output $y$
\begin{equation}
\label{eq:Txw}
T_{y,w}(s)=\frac{C(s)G(s)}{1+C(s)G(s)\frac{1}{D(s)}}F(s)=\frac{C(s)G(s)}{D(s)+C(s)G(s)}D(s)F(s),
\end{equation}
and from the system output disturbance $d$ to the system output $y$
\begin{equation}
\label{eq:Txw}
T_{y,d}(s)=\frac{1}{1+C(s)G(s)\frac{1}{D(s)}}F(s)=\frac{1}{D(s)+C(s)G(s)}D(s)F(s).
\end{equation}
Note that the term $D(s)F(s)$ appears in both the transfer functions, 
which confirms that the oscillatory mode compensation takes place in both of the channels.
Furthermore, note that the characteristic equation determining the closed-loop poles is in both cases given by
\begin{equation}
\label{eq:CE}
D(s)+C(s)G(s)=0.
\end{equation}
If $C(s)G(s)$ is a strictly-proper transfer function, 
then $\lim_{|s|\rightarrow \infty}|C(s)G(s)|=0$, and so the closed-loop poles with large magnitudes will tend to match the spectrum of $D(s)$ zeros. 
For more details, particularly if $C(s)G(s)$ is bi-proper, we refer the reader to \cite{vyhlidal2015feedback}.

As designing a robust input shaper can lead to the shaper having a large number of zeros in the right half-plane (see \cite[Section 2.2]{KunPVetal17}),
designing $C(s)$, the controller for the system $G(s)\frac{1}{D(s)}$, would also have the additional task of restricting all the 
poles of the closed-loop system to the open left half-plane.
However, the inclusion of the aforementioned stability constraint to restrict the zeros of the shaper to open left half-plane greatly reduces this burden on the controller.
This leads to defining the so-called \emph{zeros spectral abscissa} of the input shaper $D(s)$:
\begin{equation}
	\alpha_D \coloneqq \max \{\Re \lambda : \text{for all } \lambda \in \C \text{ such that } D(\lambda) = 0\}.
\end{equation}
In this paper, we only consider stability from the perspective of the zeros of the shaper;  
extension to the full closed-loop system has not yet been in the consideration in the literature
but will be the subject of future work.

We now present the specific optimization problem for designing the gains of the delay spikes in the input shaper.
As its derivation is both technical and not directly relevant to our goal of how best 
to solve it in practice, we omit many of these details here and instead just refer the reader to \cite{PilMV17} for how they exactly arise.
Given~\eqref{eq:shaper_dyn} as the chosen realization of the shaper, with defining constants $T$, $n$, and $\initgain$,
the optimal \emph{stable} shaper is obtained via solving
\begin{equation}
\label{qcconvex2}
\min_{x\in\R^n} \, 
	\frac{	
		x\tp H x 
	}{
		R_\mathrm{nom}^2
	}  \quad \text{s.t.} \quad 
	\begin{cases} 
	A_1 x & \leq  b_1 \\
	A_2 x & =  b_2 \\
       	\zsa(x) & \le  \zsamargin
\end{cases},
\end{equation}
where $A_1 \in \Rmn{(n-1)}{n}$ and $b_1 \in \R^n$ in the linear inequality constraints are defined by
\begin{equation}
	\label{eq:lin_ineq}
	-\sum_{k=1}^l x_k \le 0 \quad \text{for} \quad  l \in \{1,\ldots,n-1\} ,
\end{equation}
$A_2 \in \Rmn{2}{n}$ and $b_2 \in \R^2$ in the linear equality constraints are defined by
\begin{equation}
	\label{eq:lin_eq}
	\sum_{k=1}^n x_k = 0 
	\qquad \text{and} \qquad 
	\gamma - \sum_{k=1}^n x_k\tau_k = 1, 
\end{equation}
$H \in \Rmn{n}{n}$ is positive definite,
$R_\mathrm{nom}>0$, 
$\zsa(x)$ is the zeros spectral abscissa of the shaper $D(s)$ as a function of optimization variables $x$,
and constant $\zsamargin \le 0$.
Matrix $H$ is a discretized approximation of the residual vibrations in the system
while $R_\mathrm{nom}$ is the reference for residual vibrations.
The linear inequality constraints \eqref{eq:lin_ineq}
ensure the desired non-decreasing step response of the shaper while 
the linear equality constraints \eqref{eq:lin_eq}
guarantee that the static gain of the shaper $D(s)$ remains equal to one,
so as to not otherwise influence the system. 
The choice of the nonpositive constant $\zsamargin$ 
allows one to optionally enforce that the zeros of the shaper not just 
merely be in the left half-plane but that they be at least $ |\zsamargin |$ away from the imaginary axis.

\begin{rem}
We note that we present \eqref{qcconvex2} with a few technical/notational differences 
from what is used in the preceding control literature.  
First, in order to maintain consistency with \eqref{eq:qp_nsineq} 
and the convention of the solvers in our evaluation,
we flipped the $\geq$ relation 
for linear inequalities to $\leq$, i.e. by implicitly multiplying them by $-1$.
Second, we use $\initgain$ to denote the initial gain, instead of $A$ as has been done elsewhere.
Third, in previous literature, the initial gain $\initgain$ was embedded as the first element of a vector $\textbf{x}$ also containing the individual gains to be optimized,
i.e. $\textbf{x} \coloneqq \begin{bsmallmatrix} \initgain & x_1 & \ldots & x_n \end{bsmallmatrix}\tp$.
Rewriting the problem with $\initgain$ removed from $\textbf{x}$ is straightforward
and can still be rewritten as the quadratic program plus a nonsmooth constraint we give in \eqref{qcconvex2}.
However, to accommodate this (small) reformulation of the vector of gains, 
it means the matrices we refer to above are correspondingly \emph{transformed versions}
of those defined in the aforementioned control literature, and are thus not the same.
\end{rem}

\subsection{The zeros spectral abscissa}
As a max function, the zeros spectral abscissa of a time-delay system will, in general, be nonsmooth when there are ties, 
that is when there is more than a single globally rightmost zero (ignoring conjugacy).
Furthermore, as the spectral abscissa of a linear system is known to be non-Lipschitz when an eigenvalue 
attaining the spectral abscissa has multiplicity greater than one \cite{BurO01}, 
the zeros spectral abscissa is also not Lipschitz, since a linear system is simply a special case of a time-delay system.
As such, one should expect that optimizing the spectrum of a time-delay system will, in general, be at least as hard as that of a linear one.

To compute the zeros spectral abscissa of $D(s)$ as a function of the optimization variables $x$, that is $\zsa(x)$,
we actually instead consider computing roots of $sD(s)$, i.e.
$\lambda \in \C$ such that
\begin{equation}
	\label{eq:shaper_spectrum}
 	\lambda D(\lambda) = \lambda\initgain + \sum_{k=1}^n x_k e^{-\lambda \tau_k} = 0.
\end{equation}
The reason for this is that the roots of equation \eqref{eq:shaper_spectrum}, which is of retarded quasi-polynomial form, can be determined as the poles of the associated retarded time-delay system
\begin{equation}
	\label{eq:assocsystem}
	\dot{y}(t)=-\frac{1}{\gamma}\sum_{k=1}^n x_k y(t-\tau_k),
\end{equation}
and so the roots of \eqref{eq:shaper_spectrum} can be computed via the method of \cite{wu2012reliably}, which is available as an open-source MATLAB code.\footnote{\url{http://twr.cs.kuleuven.be/research/software/delay-control/roots/}}
Note that considering $sD(s) = 0$ introduces an additional root at the origin, 
but this is simply removed from the computed set of roots.

In general, \eqref{eq:shaper_spectrum} has infinitely many solutions. However, due to $\gamma>0$, the spectrum is retarded. 
For retarded systems, 
the number of roots satisfying $\Re \lambda >\beta$ for any $\beta\in\mathbb{R}$ is always finite \cite{hale2013introduction}. 
This is in constrast to neutral time-delay systems, 
where the number of roots located to the right of a given vertical line specified by $\beta$ can be infinite. 
Thus, a neutral system can have infinitely many unstable roots,  
which in turn makes stabilizing them more difficult than 
retarded systems, which only have finitely many unstable roots.
Note that a classical shaper with lumped delays has a neutral spectrum. 
Avoiding this undesirable spectrum distribution, especially for the case when the shaper is applied in the inverse form, 
was the main reason why distributed-delay shapers were introduced \cite{vyhlidal2013signal}.

Although apparently not yet rigorously proven, 
the zeros spectral abscissa at least seems to be differentiable almost everywhere with respect to the delay gains $x$,
even though it is still ultimately nonsmooth \cite[Section 10.2.1]{MicN07}.
Let $\hat x$ be any set of specific gains such that $\zsa(x)$ is smooth at $\hat x$, 
with $\lambda_1 \in \C$ being a simple globally rightmost zero of \eqref{eq:shaper_spectrum} (again, excluding conjugacy).
It then follows from straightforward differentiation with respect to $x$ that the gradient of the zeros spectral abscissa at $\hat x$ is given by
\begin{equation}
\label{eq:spabs_grad}
\nabla \zsa(x) \bigg\vert_{x = \hat x} = - \Re\left( 
\frac{
\begin{bmatrix}
	e^{-\lambda_1 \tau_1} & \ldots &  e^{-\lambda_1 \tau_n}
\end{bmatrix}\tp
}
{
\initgain - \sum_{k=1}^n\tau_k \hat{x}_k e^{-\lambda_1 \tau_k} 
}
\right).
\end{equation}
Although computing the zeros of the shaper can be nontrivial,
the subsequent gradient calculation is relatively straightforward and inexpensive to obtain.

\section{Descriptions of SQP-GS and GRANSO}
\label{sec:solvers}
For brevity in describing the algorithms, 
we will limit the discussion in this section to inequality constrained optimization 
but we note that both methods can also handle equality constraints and 
indeed, for the optimization problems we wish to solve here, 
we require methods that are applicable for problems with both equality and inequality constraints.
Assuming that all the functions in \eqref{eq:nlp} are locally Lipschitz, 
by Rademacher's Theorem~\cite{RocW09},
it follows that any nondifferentiability must be restricted to sets of measure zero, if at all.
In our setting here, 
we will assume that there is always at least one function in \eqref{eq:nlp} that is nonsmooth. 
Nevertheless, even when the nonsmoothness of \eqref{eq:nlp} is restricted to just a set of measure zero,
this is still often an insurmountable challenge to overcome for many optimization methods 
that are not specifically built for nonsmooth optimization in mind.
In particular, in many applications, 
minimizers of \eqref{eq:nlp} will typically be points at which one of the functions is nondifferentiable.

One potential way to solve \eqref{eq:nlp} is by minimizing a (nonsmooth) \emph{penalty function},
an approach employed by both SQP-GS and GRANSO.
Consider the $l_1$ penalty function 
\begin{equation}
	\label{eq:penfun}
	\phi(x; \rho) \coloneqq \rho f(x) + v(x)
\end{equation}
where penalty parameter $\rho > 0$ balances the (sometimes opposing) goals 
of minimizing the objective and satisfying the constraints
and $v : \R^n \rightarrow \R$ measures 
the total violation of constraints, defined as
\begin{equation}
	\label{eq:violfun}
	v(x) \coloneqq  \| \max\{c(x),0\} \|_1  =  \sum_{j \in \mathcal{P}_x} c_j(x) 
		\quad \text{where} \quad 
	\mathcal{P}_x \coloneqq \{j \in \{1, \dots, m\} \text{ : } c_j(x) > 0 \}.
\end{equation}
If an acceptable value of $\rho$ is known \emph{a priori},
then solving \eqref{eq:penfun} via a suitable \emph{unconstrained} optimization method
will also yield a solution to the constrained problem \eqref{eq:nlp} (see \cite{Bur91,Bur91a}).
Of course in practice, this is often not the case and the responsibility 
of finding a good value for $\rho$ falls upon the constrained optimization method.

\subsection{SQP-GS}\label{sec:sqpgs}
In order to overcome any discontinuities in the gradients 
of the objective or constraint functions, respectively $\nabla f$ and $\nabla c^j$,
SQP-GS finds a stationary point of \eqref{eq:nlp} by solving a sequence of quadratic programs
augmented by additional gradient information obtained by random sampling.
At each iteration,
$\nabla f$ and $\nabla c^j$ are randomly sampled from $\epsilon$-neighborhoods around the current iterate $x_k$; 
we denote these respective sets of sampled gradients by 
$\mathcal{B}^f_{\epsilon,k}$ and $\{\mathcal{B}^{c^j}_{\epsilon,k}\}_{j=1}^m$.
For now, we will assume that the objective and constraint functions are all nonsmooth.
In order for the convergence results of SQP-GS to hold, 
each set must contain at least $n+1$ sampled gradients (in addition to the gradient at $x_k$) but in practice, 
it is usually recommended to sample more than this (e.g. $2n$).
Then to compute a descent direction $d$, the following subproblem of \eqref{eq:nlp} is solved:
\begin{equation}
\label{eq:sqpgs_dir}
	\begin{aligned}
		\min_{d,z,r} &\quad \rho z+\sum_{j=1}^m r^j+\frac{1}{2}d\tp H_k d \\
		\text{s.t.} &\quad f(x_k)+\nabla f(x)\tp d\le z,\ \ \forall x\in \mathcal{B}^f_{\epsilon,k} \\
			&\quad c^j(x_k)+\nabla c^j(x)\tp d\le r^j, \ \ 
				\forall x\in \mathcal{B}^{c^j}_{\epsilon,k}, \ \ r_j \ge 0, \ \ \forall j\in\{1,\ldots,m\},
	\end{aligned}
\end{equation}
where $\rho$ is the aforementioned penalty parameter, 
$z \in R$ and $r^j \in \R$ are slack variables,
and $H_k$ is a BFGS approximation of the Lagrangian Hessian
that is maintained and updated at every iteration.
Then given the descent direction $d$ obtained by solving \eqref{eq:nlp}, 
a backtracking line search is performed to compute the next iterate $x_{k+1}$ 
such that sufficient decrease of the penalty function \eqref{eq:penfun} is satisfied. 

Assuming that all the functions in \eqref{eq:nlp} are indeed locally Lipschitz, 
then with probability one,
SQP-GS is guaranteed to converge to a stationary point of \eqref{eq:nlp}
but this computed stationary point may be either a feasible or infeasible one.
Note that the convergence results for SQP-GS also require that 
the BFGS Hessian approximation remain bounded,
which is enforced in the 1.2 and later versions of its implementation. 
As mentioned earlier, the main downside to using SQP-GS is
the many gradient samples ($\bigO(n)$) that need to be evaluated on each iteration,
which can make SQP-GS a very expensive method to use, perhaps prohibitively so,
when the functions and/or gradients are expensive to evaluate.
Unfortunately, this is often the case for applications in robust control,
including the problem we consider in this paper.
Nevertheless, it is only necessary to apply the gradient sampling procedure 
to the specific functions that are actually nonsmooth and so SQP-GS allows 
a user to set the number of samples individually for each individual function.
While we will make use of this feature to only enable gradient sampling for the nonsmooth zeros spectral abscissa,
the relative cost of computing the remaining smooth functions and their gradients
is so negligible that such selected gradient sampling will only marginally improve the computational burden for our specific problem.

\subsection{GRANSO}\label{sec:granso}
To avoid the potentially high overhead that SQP-GS can have, 
GRANSO completely eschews gradient sampling.
Instead, inspired by the wealth of evidence that regular BFGS is a reliable and efficient method
for unconstrained nonsmooth optimization (e.g. \cite{LewO13}),
GRANSO uses BFGS updating as a tool for solving nonsmooth constrained problems, in the form of \eqref{eq:nlp},
via finding a minimizer of the penalty function \eqref{eq:penfun}.
However, the choice of determining a good value for the penalty parameter $\rho$ still remains an issue.  
To address this, GRANSO combines BFGS updating with an SQP-based \emph{steering} strategy.
Originally proposed by \cite{ByrNW08,ByrLN12} for smooth problems,
steering attempts to dynamically determine at every iteration whether or not 
the current value of $\rho$ should be lowered, and by how much, 
in order to promote sufficient progress towards the feasible region.
This strategy works as follows.

Consider the following quadratic program:
\begin{equation}
	\label{eq:granso_dir}
  	\begin{aligned}
    		\min_{d,r} &\quad \rho(f(x_k) + \nabla f(x_k)\tp d) + \sum_{j=1}^m r^j + \thalf d\tp H_kd \\
    		\text{s.t.}  	&\quad c^j(x_k) + \nabla c^j(x_k)\tp d \le r^j,\ \ r^j \ge 0 , \ \ \forall j\in\{1,\ldots,m\}.
	\end{aligned}
\end{equation}
At each iterate $x_k$, in order to determine a new descent direction $d_k$, 
GRANSO will solve \eqref{eq:granso_dir} one or more times,
for different trial values of the penalty parameter.
Note if $\rho=1$ and there are no constraints, 
solving \eqref{eq:granso_dir} would just yield the standard BFGS direction.
First, GRANSO solves \eqref{eq:granso_dir} for the current value of the penalty parameter
to produce a candidate direction $d$.  
Using a linearized model of constraint violation, 
GRANSO computes the predicted reduction in total violation if the candidate
direction $d$ is used, that is,
\begin{equation}
	\label{eq:violred}
	l_\delta(d;x_k) \coloneqq v(x_k) -  \|\max\{c(x_k) + \nabla c(x_k)\tp d,0\}\|_1.
\end{equation}
If the prediction reduction given by \eqref{eq:violred} is deemed sufficiently large, 
then the computed direction $d$ and the current value of the penalty parameter
are considered to sufficiently promote progress towards the feasible region 
and $d$ is accepted as the new search direction $d_k$ while $\rho$ remains unaltered. 
Otherwise, \eqref{eq:granso_dir} is re-solved with $\rho$ temporarily set to zero
to obtain $\tilde d_k$, which is called the \emph{reference direction}.
With $\rho$ temporarily set at zero, 
the objective function and its gradient are removed from \eqref{eq:granso_dir}.
Thus, the reference direction is biased towards reducing the current violation $v(x_k)$
as much as possible regardless of what effect it may have on the objective function.
Using the predicted violation reduction for the reference direction, $l_\delta(\tilde d_k;x_k)$,
which models the best case reduction in violation one might achieve,
GRANSO then re-solves \eqref{eq:granso_dir} for successively lower values of $\rho$
until it produces a direction $d$ such that 
$l_\delta(d;x_k)$ is some sufficiently large fraction of $l_\delta(\tilde d_k;x_k)$.
In other words, in this latter scenario, GRANSO computes a direction $d_k$ that is predicted 
to reduce the violation by at least, say 10\% for example, 
of the amount that the reference direction $\tilde d_k$ would. 
An inexact Armijo-Wolfe line search is then used with direction $d_k$ to compute the next iterate
and the BFGS Hessian approximation $H_k$ is updated.

\section{Experimental setup and preliminary analysis}
\label{sec:experiments}
For setting up the specific instance of \eqref{qcconvex2}, we chose to optimize an input shaper with 18 equally-spaced delays ($n=18$) with $T=0.8$, $\initgain = 0.01$, and $\zsamargin = -0.1$ to keep the zeros of the shaper a bit away from the stability boundary.
To assess the difficulty of solving this 18-variable robust control design task and benchmark SQP-GS and GRANSO,
we ran both solvers on two different sets of 1000 starting points, each constructed as follows.
For the first set, we simply created random starting points generated via \texttt{randn},
a common (and often default) choice for most solvers and practitioners when no other information is known about the problem's solution(s) a priori.
The second set also contained 1000 randomly-generated starting points,
but created using uniform distributions such that these initial points at least satisfied all of the linear inequality constraints
and the first of the two linear equality constraints.
We will herein respectively refer to these two sets as the ``\texttt{randn}" and ``LC" points, with the latter standing for ``Linear Constraints".
We did not attempt to construct initial points that satisfied all linear constraints and/or the nonlinear zeros spectral abscissa constraint.

For the implementation of the optimization problem, 
we were provided the MATLAB code used for the experiments done in \cite{KunPVetal17} and other related papers.
However, we modified this code so that the implemented objective and constraint functions agreed 
with the formal shaper optimization problem we present here.
In this process, we also increased the code's efficiency and reliability when evaluating these function/gradient and did code cleanup and 
modularization, to facilitate our larger evaluation.

For the solvers, we specifically used SQP-GS v1.3 and GRANSO v1.6 and allowed each solver a maximum of 1000 iterations per starting point.  
We set both codes to use a stationarity tolerance of $10^{-8}$.
In order for an iterate to be considered (sufficiently) feasible, 
we set the tolerances of both solvers so that the inequality constraints had to be strictly satisfied while the total violation of the equality constraints could be no more than $10^{-8}$.
For additional consistency, 
we set the initial value of the penalty parameter of both SQP-GS (\texttt{rho\_init}) and GRANSO (\texttt{opts.mu0}) to 1.
Finally, for SQP-GS, we set the number of gradient samples for the nonsmooth spectral abscissa constraint to be $2n$,
the recommended value noted in its documentation and relevant papers (\cite[Section~4]{CurO12} and \cite[Section~4]{BurLO05}).
Otherwise, all other parameters of the two codes were left at their defaults and both codes used the builtin \texttt{quadprog} routine 
provided in MATLAB to solve their respective QP subproblems.

The experiments were computed over a subset of identical compute nodes of the \texttt{mechthild} cluster located at MPI Magdeburg.
where each node has two Intel Xeon Silver 4110 CPUs with 192 GB RAM and were running MATLAB R2017b
on CentOS Linux 7.
Note however that very little memory was actually needed for the experiments.
As both SQP-GS and GRANSO are serial programs, 
the only exploited parallelism was in the experimental setup, 
to distribute the 2000 different starting points over multiple nodes.
The cumulative serial time of all our experiments (measured using \texttt{tic} and \texttt{toc}) 
was 118 days, with SQP-GS being accountable for 97\% of the total computation time.

\subsection{Preliminary performance analysis}
In Table~\ref{tab:stats}, we provide basic cost statistics for performing the experiments, for each pair of solver and starting point set.  
The relative high cost of running SQP-GS is clearly evident, with it requiring about 31 times longer than GRANSO to complete the LC test set and needing 39 times longer on the \texttt{randn} test set.  While GRANSO completed both test sets in 3.4 days, SQP-GS took over 115 days.  
Interestingly, part of this large difference was not just the additional cost for SQP-GS to obtain $\bigO(n)$ gradient samples every iteration but also
because SQP-GS also took over twice the number of iterations compared to GRANSO.
Compared to the \texttt{randn} starting points, we see that using LC starting points reduces the total time and iterations 
for both solvers, by a factor of two, roughly speaking. 
We also see that more often than not both solvers are terminating well before their 1000 maximum allowed iteration limit. 

\begin{table}[t]
\center
\begin{tabular}{l | rr | r | rr | r}
\toprule
\multicolumn{1}{c}{} & \multicolumn{3}{c}{\texttt{randn}} & \multicolumn{3}{c}{LC} \\
\cmidrule(lr){2-4}
\cmidrule(lr){5-7}
\multicolumn{1}{c}{} & \multicolumn{1}{c}{GRANSO} & \multicolumn{1}{c}{SQP-GS} & \multicolumn{1}{c}{Ratio} & \multicolumn{1}{c}{GRANSO} & \multicolumn{1}{c}{SQP-GS} & \multicolumn{1}{c}{Ratio}\\
\midrule
Total Wall-Clock Time (days)		&	2.12 		& 83.1797 	& 39.16 	& 	1.28 		& 31.96		& 31.10	\\
Total Iters						&	107,491	& 248,254		& 2.31	&	54,151	& 118,565		& 2.19	\\
Average Time Per Iter (min/iter)		&	0.03 		& 0.48	 	& 16.96	& 	0.03		& 0.39		& 14.20	\\
\bottomrule
\end{tabular}
\caption{Total cost statistics for each solver, on each of the two test sets of 1000 different starting points.
The ``Ratio" columns simply show the statistic for SQP-GS divided by the corresponding one for GRANSO.}
\label{tab:stats}
\end{table}

\begin{figure}[t]
\center
\subfloat[Initialized from the \texttt{randn} points.]{
	\includegraphics[scale=0.43]{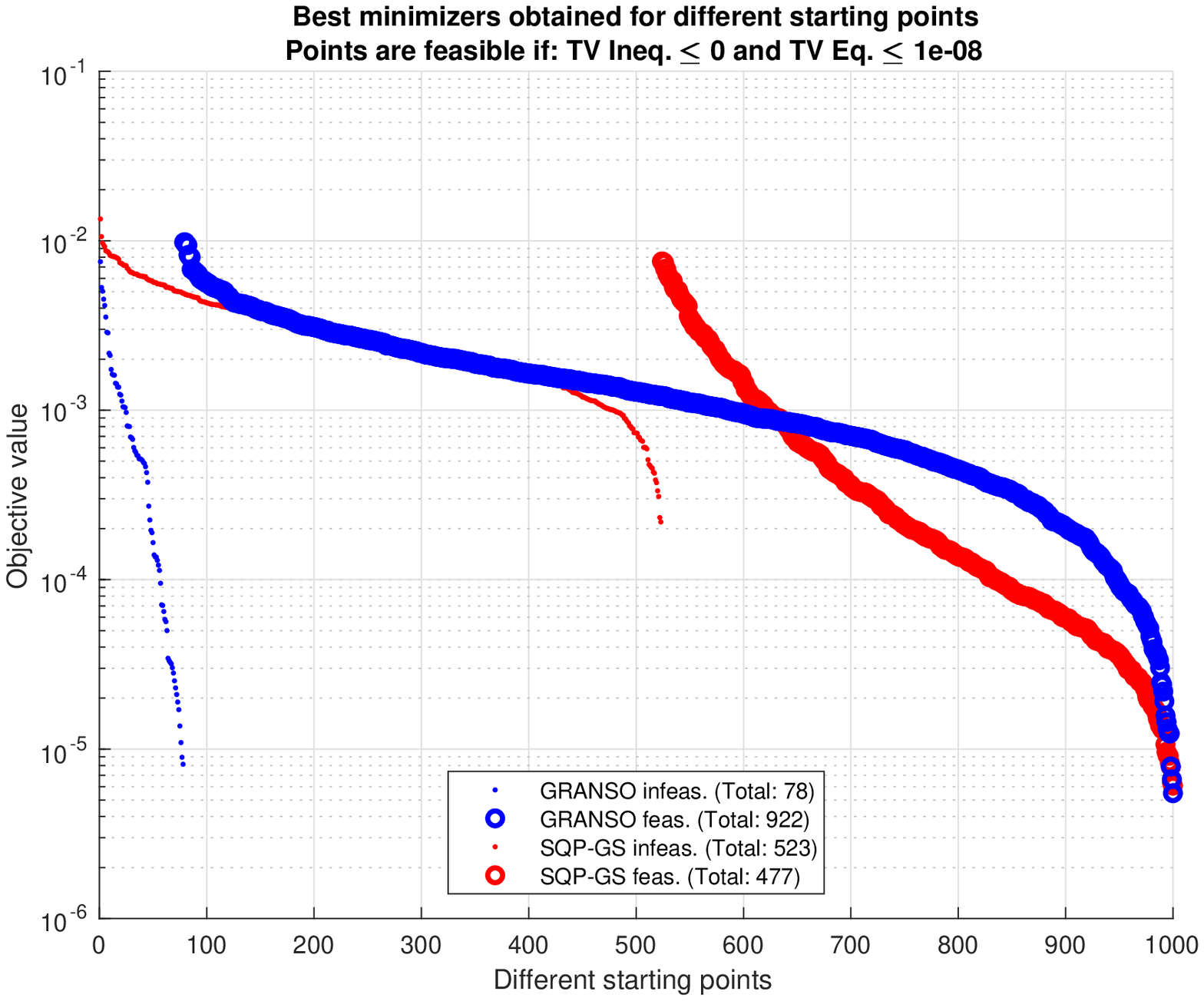}
	\label{fig-randbestvals}
}\qquad
\subfloat[Initialized from the LC points.]{
	\includegraphics[scale=0.43]{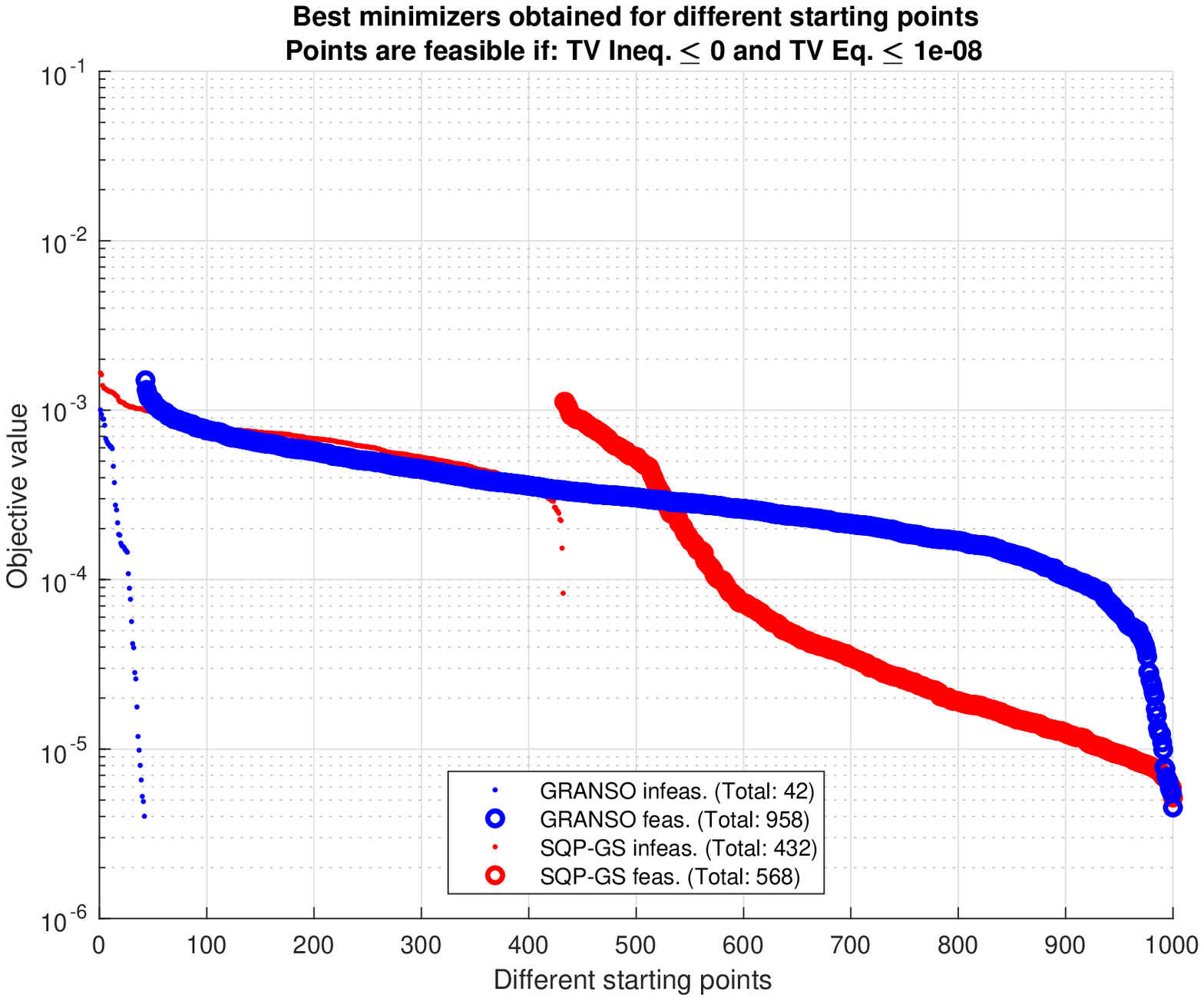}
	\label{fig-ineqbestvals}
} 
\caption{Best solutions, in terms of the value of the objective function evaluated at those solutions, obtained from each starting point.  
A solution is considered feasible if its total violation for the inequality constraints (TV Ineq.) is less than zero and its total violation for the equality constraints (TV Eq.)
is at most $10^{-8}$.
Note that use of log scale on the $y$-axis, as the best objective values span a great range.
}
\end{figure}

In Figure~\ref{fig-randbestvals}, for each of the \texttt{randn} starting points, we plot the best objective value encountered by each solver.
The plot is sorted in decreasing order with respect to these best objective values, with the infeasible solutions 
plotted first.  An infeasible solution means that a solver was not able to ever find the feasible set when initialized at that
particular starting point.  
A dramatic performance difference between GRANSO and SQP-GS is immediately apparent, namely that 
SQP-GS only finds the feasible sets for just less than half (477 of 1000) of the \texttt{randn} test set while GRANSO
finds the feasible set for 922 of the 1000 \texttt{randn} starting points.
On the other hand, 
of the feasible solutions found, we see that SQP-GS found a larger portion of the best solutions,
even though the very best solution was computed by GRANSO.

In Figure~\ref{fig-ineqbestvals}, we provide the same type of plot but now for the LC starting points.
We more or less see the same overall picture, with a couple notable differences.  First,
initializing from the LC points allows both solvers to slightly increase their respective success rates in finding the feasible region, 
which is not so surprising given the construction of the LC points to partially satisfy feasibility.  
Second, both solvers on average seem to return better solutions, ones corresponding to somewhat lower objective values,
when initialized from the LC points,
with this positive change being slightly more prominent for SQP-GS.
Nevertheless, the very best solution was also computed by GRANSO on the LC test set.

Plots like the ones we show in Figures~\ref{fig-randbestvals} and \ref{fig-ineqbestvals} are useful 
for assessing how many different stationary points a nonconvex problem may have, which would appear
as plateaus due to the sorting of the solutions by objective value.
However, neither plot seems to show any evidence of plateaus.  
Instead, both show a rather smoothly-varying decrease in the objective values.
One possible explanation is that there is only one feasible minimizer and perhaps a second infeasible one.
In this case, the smooth changes seen in the objective value could be due to inconsistencies in the computed solutions,
due to scaling issues, choices of stopping tolerances for the solvers, 
\texttt{quadprog} and its parameters, inaccuracy in the computation of the zeros spectral value abscissa, 
high sensitivity of the functions to the variables near minimizers, etc.
On the other hand, given that the zeros spectral value abscissa is likely not Lipschitz, and thus neither solver would have convergence guarantees,
the smoothly-varying objective value plots may simply be an indicator that the solvers are simply terminating before actually reaching stationary points.
In this case, it may be that there are multiple stationary points present (even though they were not detected in Figures~\ref{fig-randbestvals} and \ref{fig-ineqbestvals}). 
Unfortunately, it is difficult to determine which of the many possible combinations of explanations is correct.

\section{A comparison using Relative Minimization Profiles}
\label{sec:rmps}
The preceding analysis comparing SQP-GS and GRANSO has so far only considered the 
quality of the best answers they produced.  
Neither Figure~\ref{fig-randbestvals} or \ref{fig-ineqbestvals} takes into account 
how much computation was actually necessary to obtain these best solutions.
For assessing the relative overall efficiency and speed of progress of these two methods,
when computation is not an unlimited resource,
we will employ what are called \emph{$\beta$-relative minimization profiles} ($\beta$-RMPs or RMPs), 
new benchmarking visualization tools
first introduced in \cite[Section 7.5]{Mit14} and 
then more formally presented in \cite[Section 5]{CurMO17}.
To aid the reader who may not be familiar with these tools, 
we will first provide the prerequisite description of $\beta$-RMPs;
for brevity, we only provide a high-level overview sufficient for understanding 
our comparisons done here and 
refer to \cite[Section 5]{CurMO17} for the formal mathematical description of $\beta$-RMPs
and further notes on their usage.
Following that, we will then describe how 
we tailor the RMP benchmark for the purposes of our comparison 
(since RMPs can be considered more of a methodology, not just a single fixed specification)
and then present our numerical evaluation of SQP-GS and GRANSO.

\subsection{An RMP primer}
RMPs are typically presented as a small selection of related plots, say four, 
called a \emph{$\beta$-RMP benchmark panel},
that simultaneously compares methods in terms of amount of minimization achieved, 
ability to find feasible solutions, and speed of progress.  
Each particular RMP plot in the selection is defined for its own value of parameter $\beta$, 
which must be positive and 
scales the computational budget allowed per each problem in a given test set.
Individual budgets may be assigned to each problem differently (for example, to account 
for its relative difficulty, size, etc).
Then, these budget assignments are scaled by the choices of $\beta$ 
to determine the actual budgets for a specific profile;
the set of such RMPs are referred to as $\beta$-RMPs.

When $\beta = \infty$, the corresponding $\infty$-RMP compares methods 
without any regard to their respective computational costs.
Each method has its own \emph{$\infty$-RMP curve} on the $\infty$-RMP plot
and for a given point $(x,y)$ on such a curve, 
$y$ is the \emph{percentage of problems solved} such that the solutions
computed by this method were both feasible (to tolerances) and 
the relative differences of the corresponding objective values to 
the \emph{best known objective values} (smallest and obtained on the feasible region) 
were at most $x$ for each problem.
On the right side of the $\infty$-RMP, typically $x=\infty$ is shown, 
which gives the percentage of problems that a method was able to 
find feasible solutions for, regardless of their corresponding objective values
(since $x=\infty$ means that there is no limit on what is an acceptably small relative difference).
The left side of the $\infty$-RMP, say for $x = 10^{-16}$, 
instead gives the percentage of problems that a method found solutions
that were both feasible and agreed to the best known objective values 
(again, over the feasible set) to essentially machine precision.
Thus, the $\infty$-RMP curve for each method shows exactly how the percentage
of its \emph{acceptably solved problems} increases as the relative difference tolerance $x$ 
is \emph{loosened} from machine precision to $\infty$.  

When $\beta$ is a finite positive value, the corresponding $\beta$-RMP compares methods 
when a per-problem computational budget, scaled by $\beta$, is forced upon them.
Thus, instead of being compared with respect to their best computed answers, 
they are compared with respect to the best answers they were able to find 
within some given per-problem budgets, scaled by $\beta$.
When $\beta = 1$, the per-problem budgets are not scaled and 
the methods are judged with respect to however the per-problem budgets are specified.
As $\beta$ is increased, 
the allowable computational budgets per problem is increased,
giving each method more computation resources 
(perhaps measured in wall-clock time, number of iterations, number of function evaluations, etc) 
for a chance to improve upon their best known answers so far.

There is no prescribed choice or methodology for 
specifying the computational budgets used in $\beta$-RMPs.
Instead, the per-problem budgets are allowed to be defined in any way the user wishes, 
as different budget choices may reflect different goals of the benchmark.
Several example choices and their possible respective usage cases were
discussed in \cite[Sections 5.2 and 5.3]{CurMO17}, 
when introducing the formal description of $\beta$-RMPs.
For example, budgets could be specified a priori by the user, 
which may be appropriate if the user has either known limited computational resources or 
good estimates for the difficulties of each problem in the test set.
Alternatively, the per-problem budgets can be constructed using data collected from the evaluation itself.
In a multi-way comparison, 
one might use the average or median of the runtimes of each method for a given problem as 
the budget for that problem, as this may be a good estimate of that problem's difficulty.  
Another viable choice is to use the runtime per problem of a single method, 
which can be appropriate when one, for example, 
is a benchmarking a new method against existing others 
(and as was done in the evaluation of GRANSO \cite[Section 6]{CurMO17}).

The best known objective values used by the RMPs, called the \emph{target values}, 
can also be determined in multiple ways.
If the problems are convex and the global minimizers are known,
then the globally minimal objective values could be supplied directly by the user.
For nonconvex problems, 
minima may not be known and the methods may converge to different stationary points,
or perhaps not to any (e.g. if the problems are not Lipschitz).
As such, RMPs set the target value for a problem 
to the best (smallest) objective value encountered at a feasible point (to tolerances),
by any of the methods, at any iterate in their entire iteration histories, for the given problem.
RMPs also give the user the choice of defining these target values by restricting 
the pool of candidate iterates to consider to be those that are computed within 
$\beta$ times the allotted budget for that problem; 
this feature is called \emph{target scaling} and alters the target values to be 
the best known objective values \emph{so far}.
Enabling target scaling typically increases the separation of the different RMP curves
for each method in each particular $\beta$-RMP plot (which can aid in the visualization) 
and judges the methods with respect to results that are known to be 
\emph{achievable within the particular budgets given for that $\beta$-RMP plot}.

To produce a $\beta$-RMP benchmark panel, 
it is necessary to collect several pieces of information.
For each method-problem pair, one must collect a history of iterates.
The history must track the objective values, constraints values, and the computational cost
(in whatever metric one has chosen to use) needed to obtain each iterate by that method.
Note that the violation values can be used in lieu of constraint values 
but then one must take care the violation values across methods
are computed in identical ways.
A second of benefit of collecting the constraint values is that it allows 
the definition of acceptably feasible to be easily changed
without having to rerun the experiments.
It is worth noting that the objective values do not necessarily need to be the 
values of the objective function of \eqref{eq:nlp}; 
using a measure of stationarity is also a perfectly fine ``objective" value to track.

\subsection{Our $\beta$-RMP benchmark}
To create our $\beta$-RMP benchmark, we recorded the objective and constraint function trajectories of the methods for
every starting point.  For the elapsed time to obtain each iterate,
we simply recorded the total running time and then used the average elapsed time per iteration to estimate 
the cumulative elapsed time, as recommended by documentation of the open-source 
\texttt{betaRMP} software for making $\beta$-RMPs \cite{betaRMP}.

\begin{figure}[t]
\centering
\subfloat[$\beta=1$]{\includegraphics[width=0.40\textwidth]{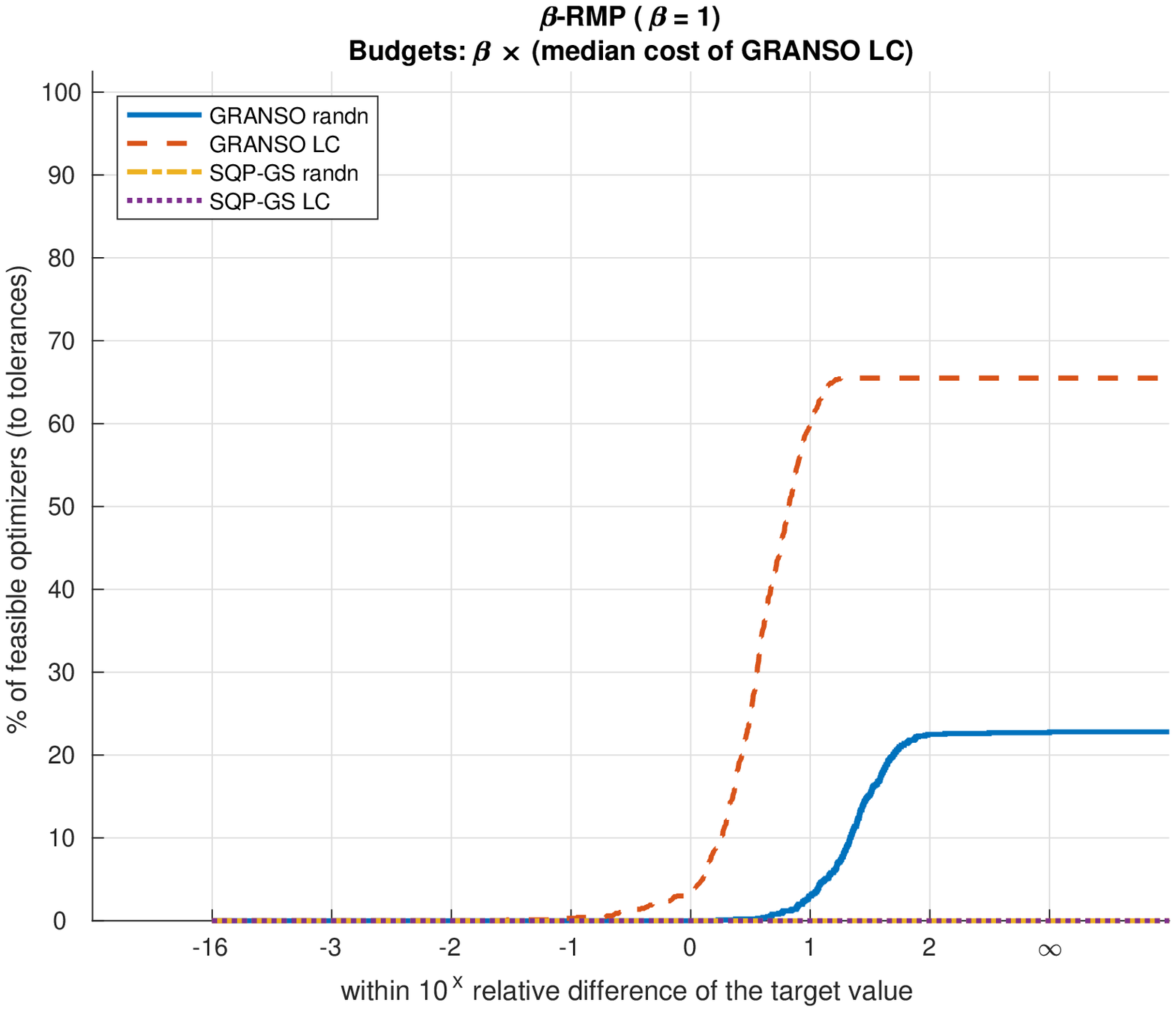}}\qquad
\subfloat[$\beta=6$]{\includegraphics[width=0.40\textwidth]{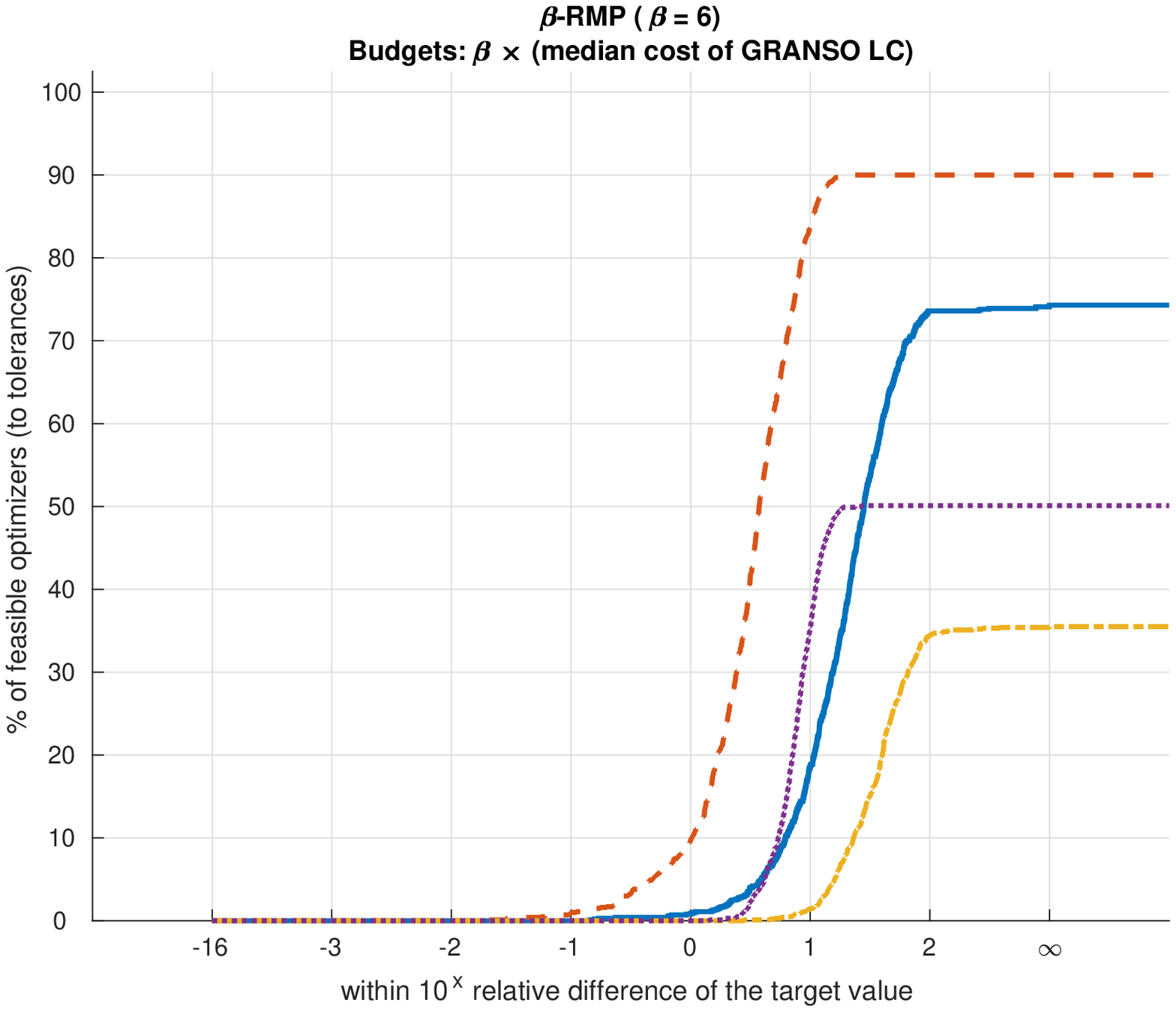}}\\
\subfloat[$\beta=12$]{\includegraphics[width=0.40\textwidth]{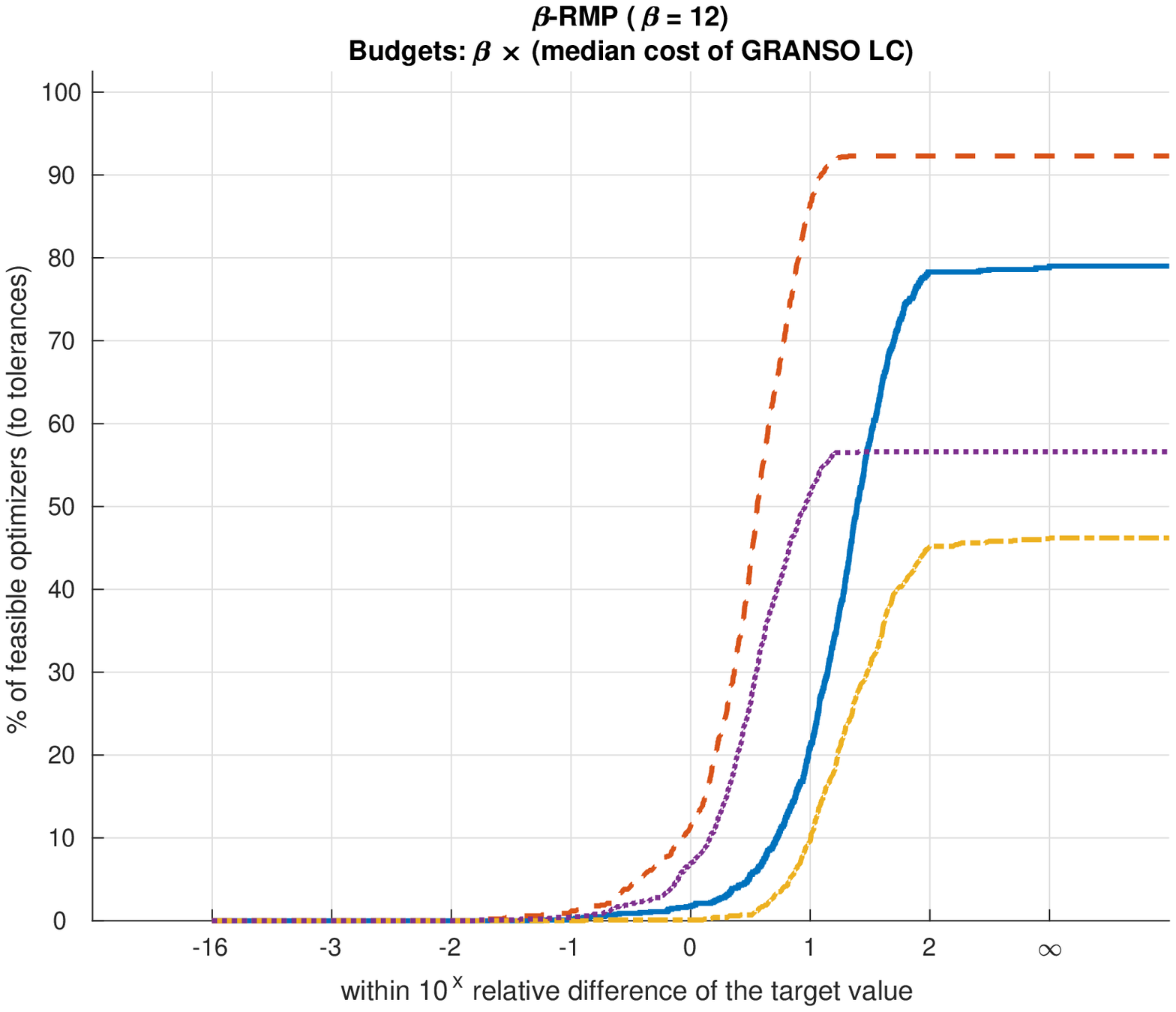}}\qquad
\subfloat[$\beta=\infty$]{\includegraphics[width=0.40\textwidth]{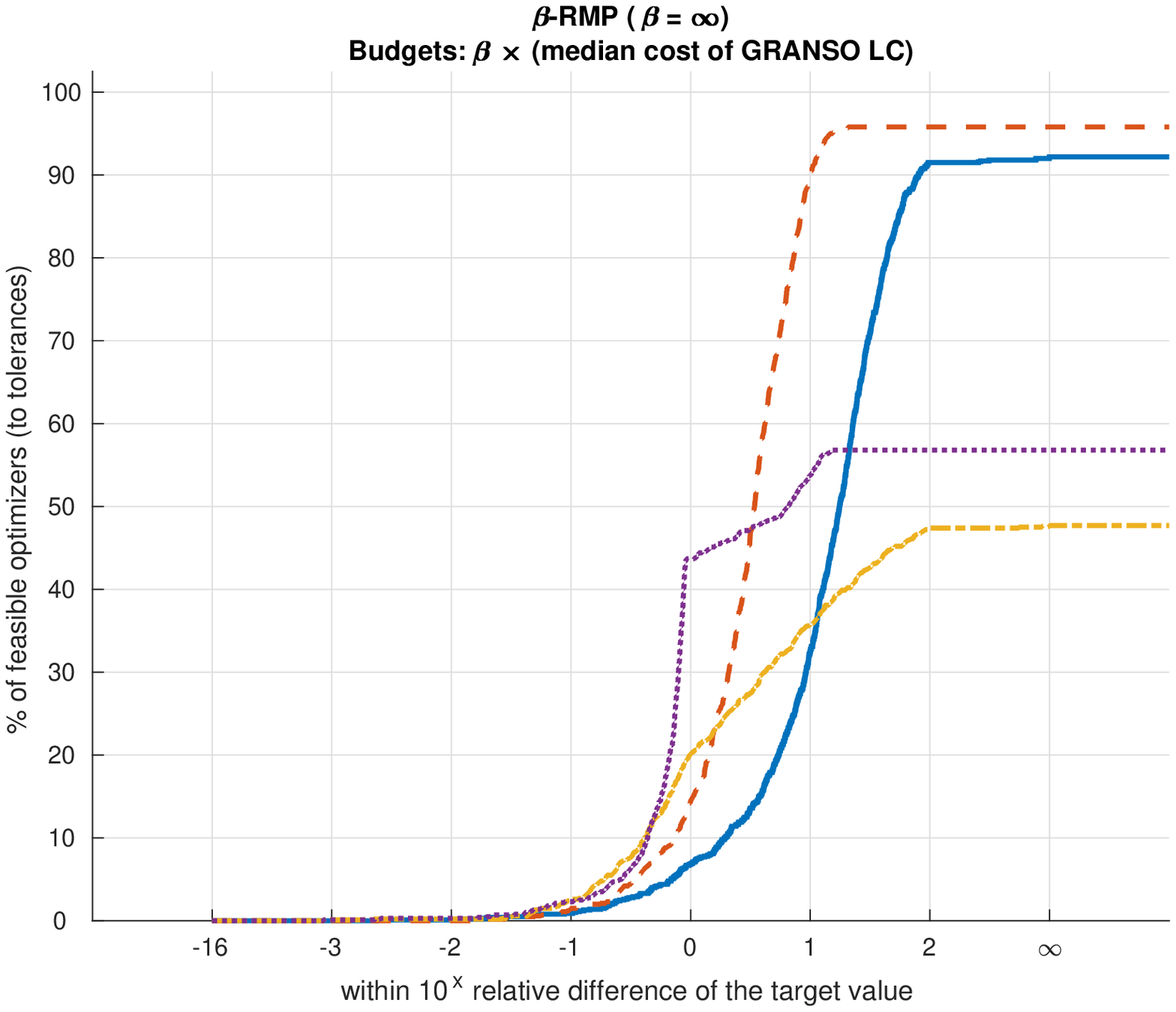}}
\caption{
A $\beta$-RMP benchmark panel for our single input shaper design problem, with a single fixed
budget and single fixed target value.
}
\label{fig:rmp2}
\end{figure}

Since we have evaluated SQP-GS and GRANSO on a single test problem, but one initialized from many different points, 
we instead specified budgets and target values a bit differently than what would be done in a RMP benchmark
over a heterogenous set of test problems.
First, we chose a single fixed budget value (to be scaled by $\beta$) for each starting point, 
since initialization technically does not change the underlying optimization problem.
For this, we set the fixed budget value as the median of the GRANSO running times over the LC points (about 25.6 seconds), 
as these GRANSO LC runs
were, on average, significantly shorter than the times of GRANSO on the \texttt{randn} points or SQP-GS on either set.
Second, we forwent the use of target scaling and simply specified a single target value,
namely the median of all the best solutions computed by both solvers, over all starting points, computed at any time;
this target value is $6.915 \times 10^{-5}$.
We selected this particular target since very few runs came close to the best computed solution, as seen in 
Figures~\ref{fig-randbestvals} and \ref{fig-ineqbestvals}.  
Third, we plotted separate RMP curves for both sets of starting points, for both methods, since there were 
noticeable differences in each of the methods' behaviors depending on which starting set was used.
Finally, we found that $\beta = \{1,6,12,\infty\}$ as the selected choice of scalings for our $\beta$-RMP benchmark panel 
best highlighted differences between the methods for various computational budgets.
The realized $\beta$-RMP benchmark panel is given in Figure~\ref{fig:rmp2}.

Remarkably, we see that even when the methods are allowed up to 12 times the median running time of GRANSO on the LC points,
SQP-GS is simply not at all competitive with GRANSO (when comparing the methods for the same initial point sets).
This is in stark contrast to the limited perspective we obtained from the best solutions plots in Figures~\ref{fig-randbestvals} and \ref{fig-ineqbestvals},
where we observed that ultimately, SQP-GS did find the majority of the best solutions if not the very best.
However, the $\beta$-RMP panel tells us that if we did not have such extended computational resources, 
SQP-GS would not have fared well at all; SQP-GS is essentially only doing well with respect to GRANSO when a concern for the cost
of the computations is marginalized.
Of course, for the smaller values of $\beta$, SQP-GS fares even worse, with the $\beta=1$ $\beta$-RMP plot showing
that even for such a limited budget, GRANSO can still frequently return sufficiently good solutions while SQP-GS returns \emph{none}.  
However, it is worth noting that for the larger values of $\beta$, we do see that SQP-GS initialized with the LC points
does start to outperform GRANSO initialized from the \texttt{randn} set, but only in terms of objective minimization and not 
in terms of success rate in finding the feasible set.

\section{Global-Local Profiles}
\label{sec:glprofiles}
For a user of optimization software, typically there is a particular desired window of time in which they would like to obtain the best possible solution 
to a problem.  For nonconvex problems, one important consideration is how much effort to devote to ``global" versus ``local" optimization,
that is, to either distribute a fixed budget over a large number of starting points or only a few.
If there are a large number of local minimizers of significantly different minimal objective values, 
then it is expected that the ideal choice might be to start from as many points as possible, with the hope that a solver
might encounter at least one of the better local minimizers, even if only approximately.
Otherwise, if there are few local minimizers 
or they all have similar minimal objective values (as is often the case in machine learning applications, for instance), 
it would likely be better to use only few starting points but to obtain these minimizers to a higher precision.

While such considerations
will partially depend on an algorithm's speed and accuracy, 
it will primarily depend on the problem itself. 
The quantity of local minimizers and the complexity and curvature of the regions around each local minimizer will vary problem to problem, and thus reward different strategies of budget distribution.
To address this problem, we propose new visualization tools called Global-Local Profiles (GL-Profiles).

In order to create GL-Profiles, we will need the same collected data as needed for $\beta$-RMPs.  
An optimization algorithm has been run $N$ times from $N$ different starting points,
with all the objective and constraint values recorded for every iterate, along with the cumulative cost to obtain each iterate;
the cumulative cost is with respect to each starting point, i.e. the cost is zero at every $x_0$.
Consider the following basic definitions regarding all the iterates of a single method started from $N$ initial points:
\begin{align*}
	x_0^{(i)} & \coloneqq \text{the $i$th starting point in the set of $N$ initial points} \\
	\{ x_k \}^{(i)} & \coloneqq \text{the sequence of iterates computed by the given solver initialized at } x_0^{(i)} \\
	t_i(j) & \coloneqq \text{cumulative cost to compute } \{x_0^{(i)},x_1,\ldots,x_j\} \subseteq \{x_k\}^{(i)}.
\end{align*}
Now, for a solver initialized at $x_0^{(i)}$, we define the subset of (sufficiently) \emph{feasible iterates} that were computed within a fixed computational limit given by $t > 0$:
\begin{equation*}
	\mathcal{F}_i \coloneqq \{ x_j : x_j \in \{x_k \}^{(i)}, \, x_j \text{ is considered feasible, and } t_i(j) \le t \}.
\end{equation*}
Let $S_l \subseteq \{1,\ldots,N\}$ denote some selection of the initial points, taken by their indices;
the subscript $l$ will come into effect later, when we consider multiple selections.  Then we consider
the subset of those for which the solver computed at least one 
feasible iterate within the computational limit given by $t$:
\begin{equation*}
	S_l^\mathcal{F} \coloneqq \{ i : i \in S_l \text{ and } \mathcal{F}_i \ne \varnothing \}.
\end{equation*}
Given some total budget $T$ and a selection of starting points given by some set of indices $S_l$, 
the budget allotted per starting point is simply $\tfrac{T}{| S_l |}$.
Then, the best (lowest) objective value obtained on the feasible set by running a solver from the points given by $S_l$, each with budget $t = \tfrac{T}{| S_l |}$, is
\begin{equation*}
	f_l^\star \coloneqq \min_{i \in S_l^\mathcal{F}} \min \{ f(x) : x \in \mathcal{F}_i \},
\end{equation*}
where $f_l^\star$ is taken to be $\infty$ if $S_l^\mathcal{F}$ is empty.

Let $P \coloneqq \{S_1,\ldots,S_q\}$ be a set of equally-sized nonoverlapping subsets of the $N$ starting points, that is,
with $| S_i | = | S_j |$ and $S_i \cap S_j = \varnothing$ for all $i \ne j \in \{1,\ldots,q\}$.
We now define the subset of sets $S_l \in P$ that have nonempty associated $S_l^\mathcal{F}$ subsets:
\begin{equation*}
	P^\mathcal{F} \coloneqq \{ l : S_l \in P \text{ such that } S_l^\mathcal{F} \ne \varnothing \}.
\end{equation*}
In other words, $P^\mathcal{F}$ is the subset of $P$ such that $f_l^\star \ne \infty$ for all $l \in P^\mathcal{F}$.
Thus, given some $P$, we can now calculate a statistic for the expected best (lowest) objective value obtained among the feasible iterates
if we were to run the solver initialized at $M = | S_1 |$ different starting points, with each given a budget of $t = \tfrac{T}{M}$:
\begin{equation*}
	\tilde f^\star \coloneqq \frac{1}{| P^\mathcal{F} |} \sum_{l \in P^\mathcal{F}} f_l^\star,
\end{equation*}
along with the standard error for this statistic:
\begin{equation*}
	\tilde f^\mathrm{err} \coloneqq \frac{1} {\sqrt{| P^\mathcal{F} |}} \sqrt{ \frac{\sum_{l \in P^\mathcal{F}} (f_l^\star - \tilde f^\star )^2} { | P^\mathcal{F} | - 1 } },
\end{equation*}
i.e. the usual standard error for estimating a mean of a value given a set of sample means, 
with estimated (rather than a priori known) population standard deviation.
Thus, given a fixed budget, a \emph{GL-Profile curve} for a method plots $\tilde f^\star$ as the number of starting points is increased,
along with the error bars given by $\tilde f^\mathrm{err}$.  
Of course, as the number of starting points is increased, the number of possible subsets in $P$ will decrease, and thus 
error bars will grow, as the $\tilde f^\star$ statistic becomes less reliable.

We can also create GL-Profile curves with error bars for just the expectation of finding the feasible set, regardless of the objective value:
\begin{equation*}
	\tilde c^\star \coloneqq \frac{| P^\mathcal{F} | }{ | P | } ,
\end{equation*}
also along with standard error for this statistic:
\begin{equation*}
	\tilde c^\mathrm{err} \coloneqq \sqrt{ \frac{\tilde c^\star (1 - \tilde c^\star)}{| P |} }.
\end{equation*}

As one might expect, a GL-Profile curve will generally be more or less ``U"-shaped.  
If few runs are performed, then the risk is high for them to be bad starting points, that either only result in infeasible iterates and/or 
high objective values.  However, this risk can quickly decrease as more initial points are added.
But, if too many points are used, any increased benefit will be outweighed by the effect of the decreasing budget allotted per starting point.
At the very extreme, with little budget devoted to each starting point, 
running optimization might not be so different than just randomly sampling the variable space.
Often, the optimal strategy, denoted by a global minimum of GL-Profile curve, will be somewhere in between the extremes but
this will depend on both the available total budget, the particular solvers, and the properties of the problem itself.

\subsection{Our GL-Profile benchmark}
\begin{figure}
\centering
\subfloat{\includegraphics[width=0.40\textwidth]{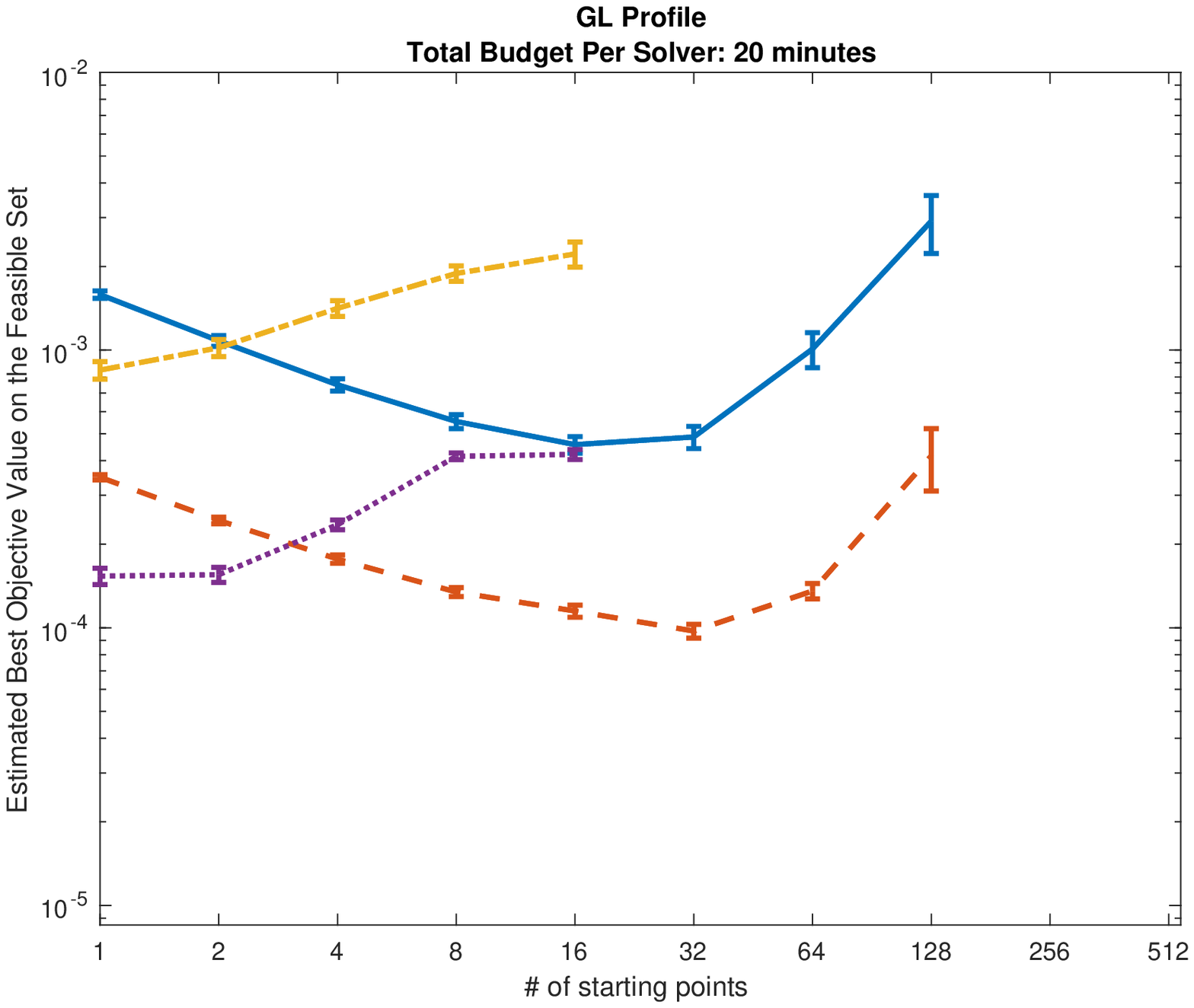}}\qquad
\subfloat{\includegraphics[width=0.40\textwidth]{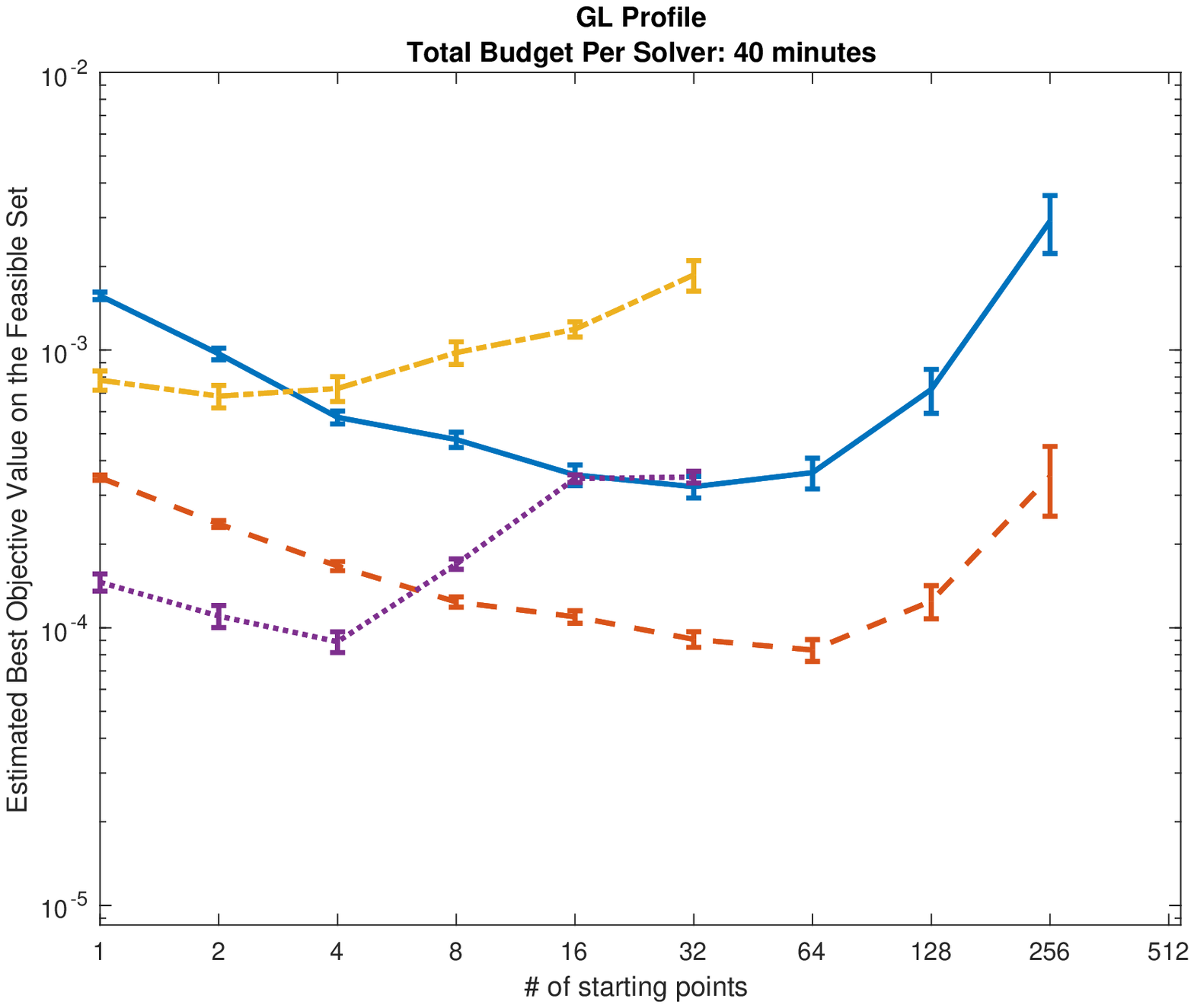}}\\
\subfloat{\includegraphics[width=0.40\textwidth]{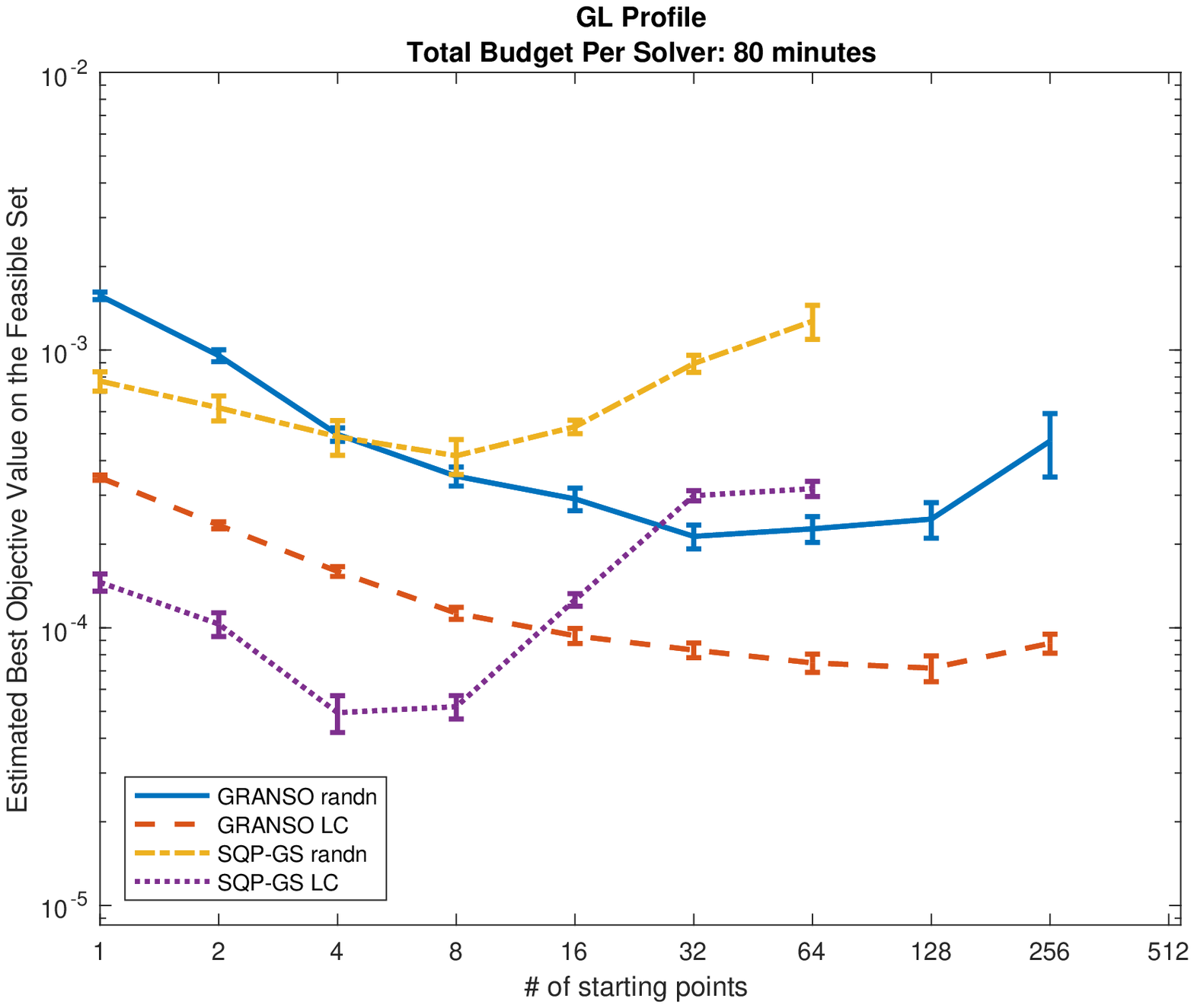}}\qquad
\subfloat{\includegraphics[width=0.40\textwidth]{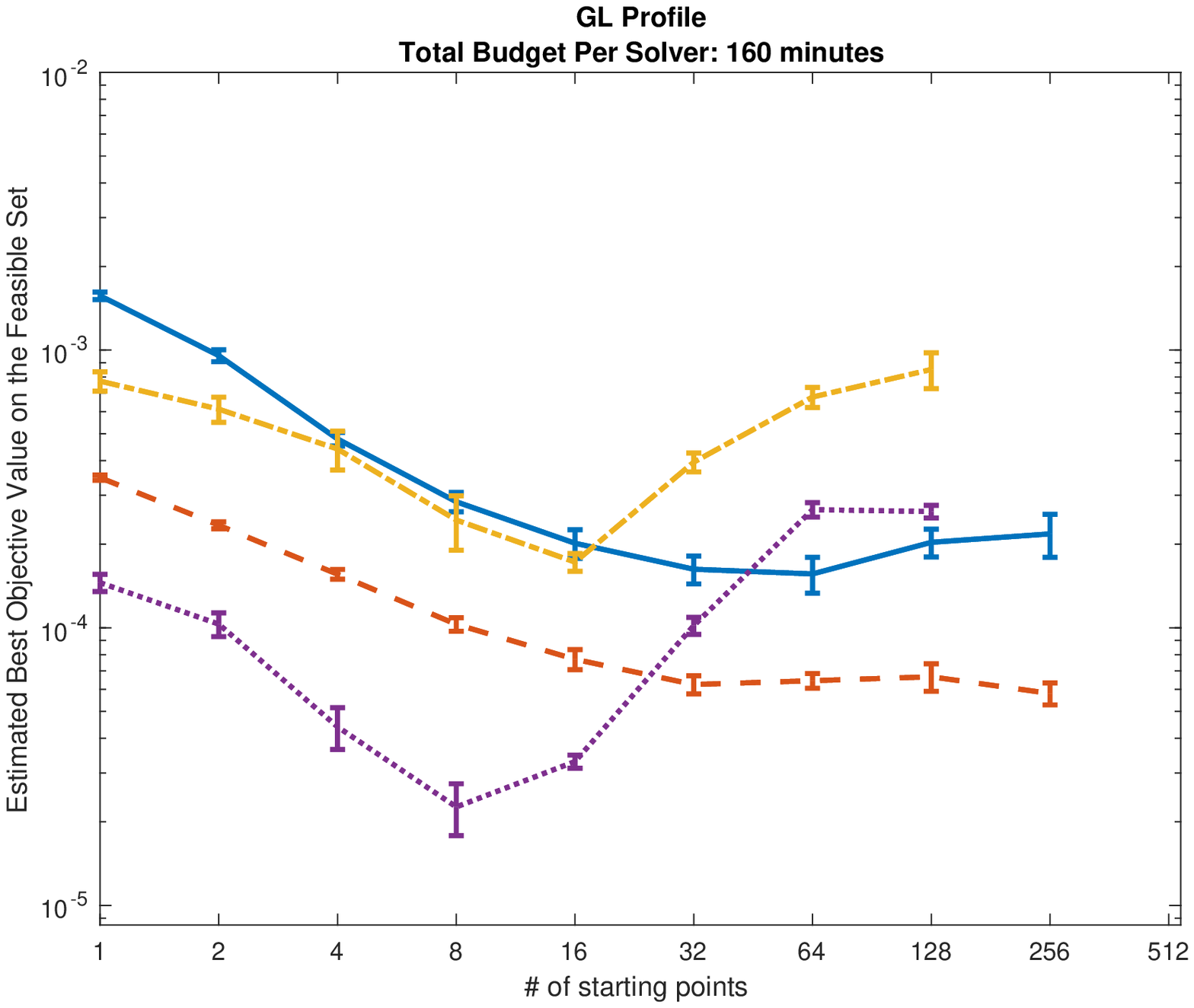}}
\caption{GL-Profiles for estimated best objective value found on the feasible set.}\label{fig:glone}
\end{figure}

\begin{figure}
\centering
\subfloat{\includegraphics[width=0.40\textwidth]{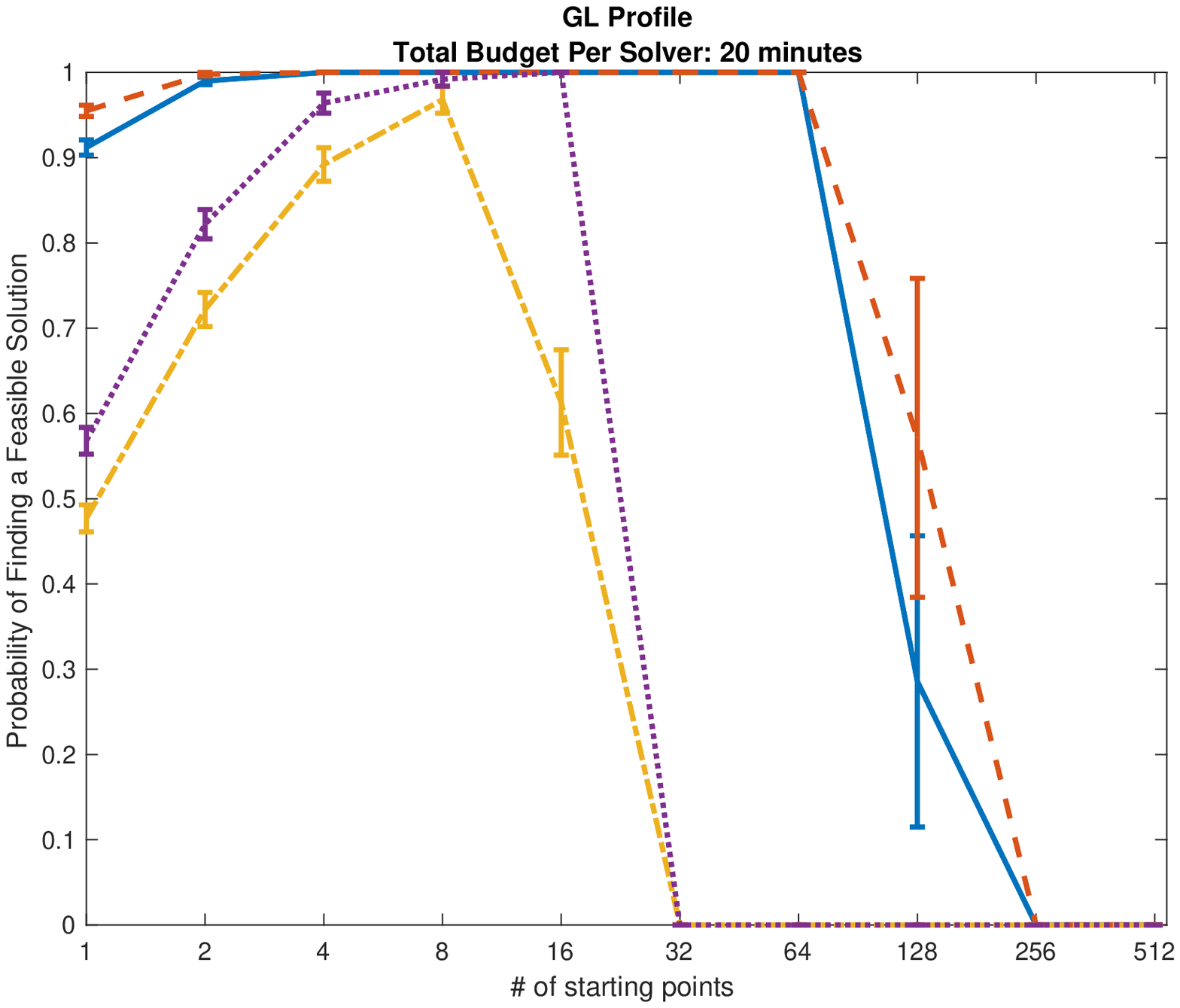}}\qquad
\subfloat{\includegraphics[width=0.40\textwidth]{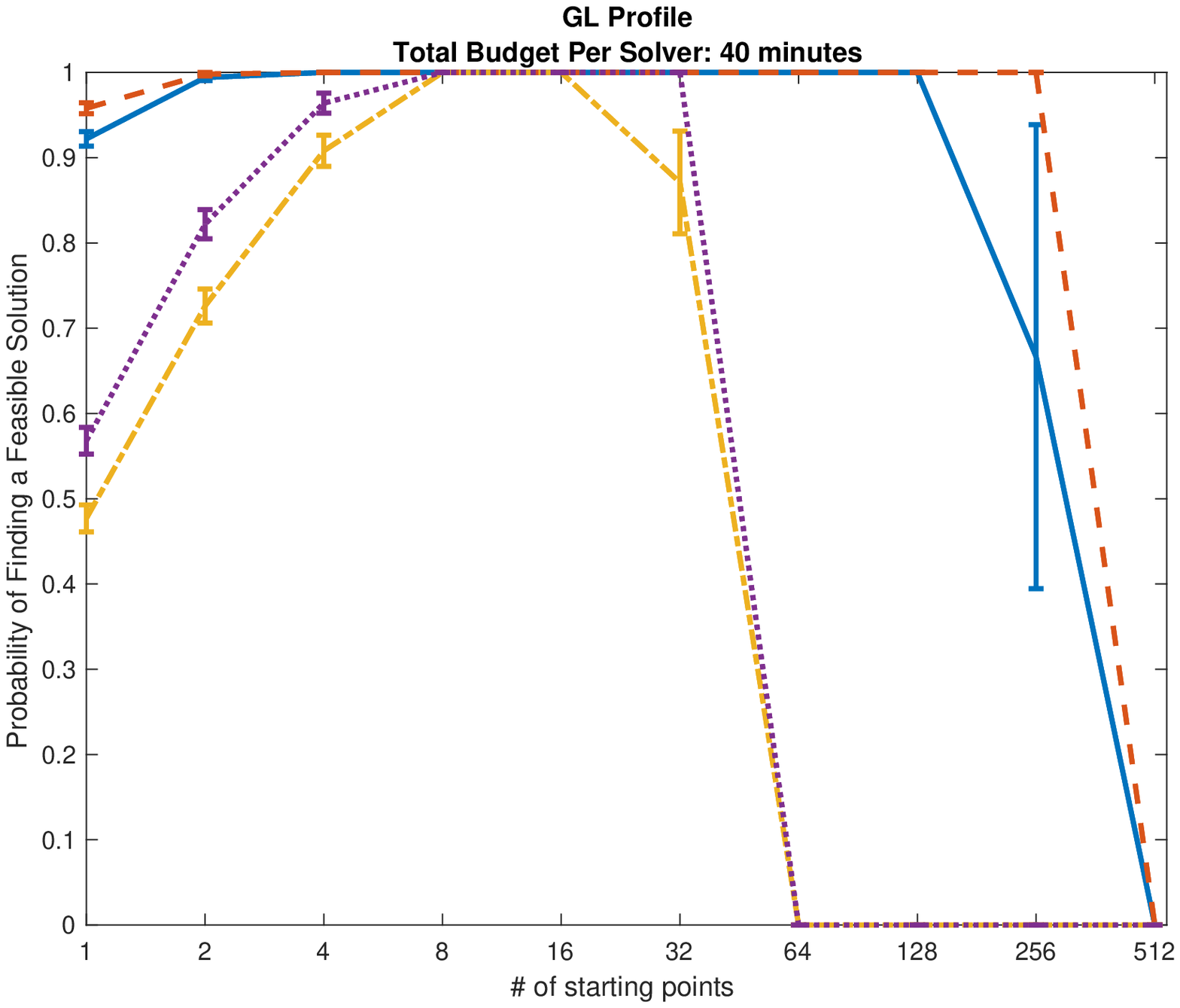}}\\
\subfloat{\includegraphics[width=0.40\textwidth]{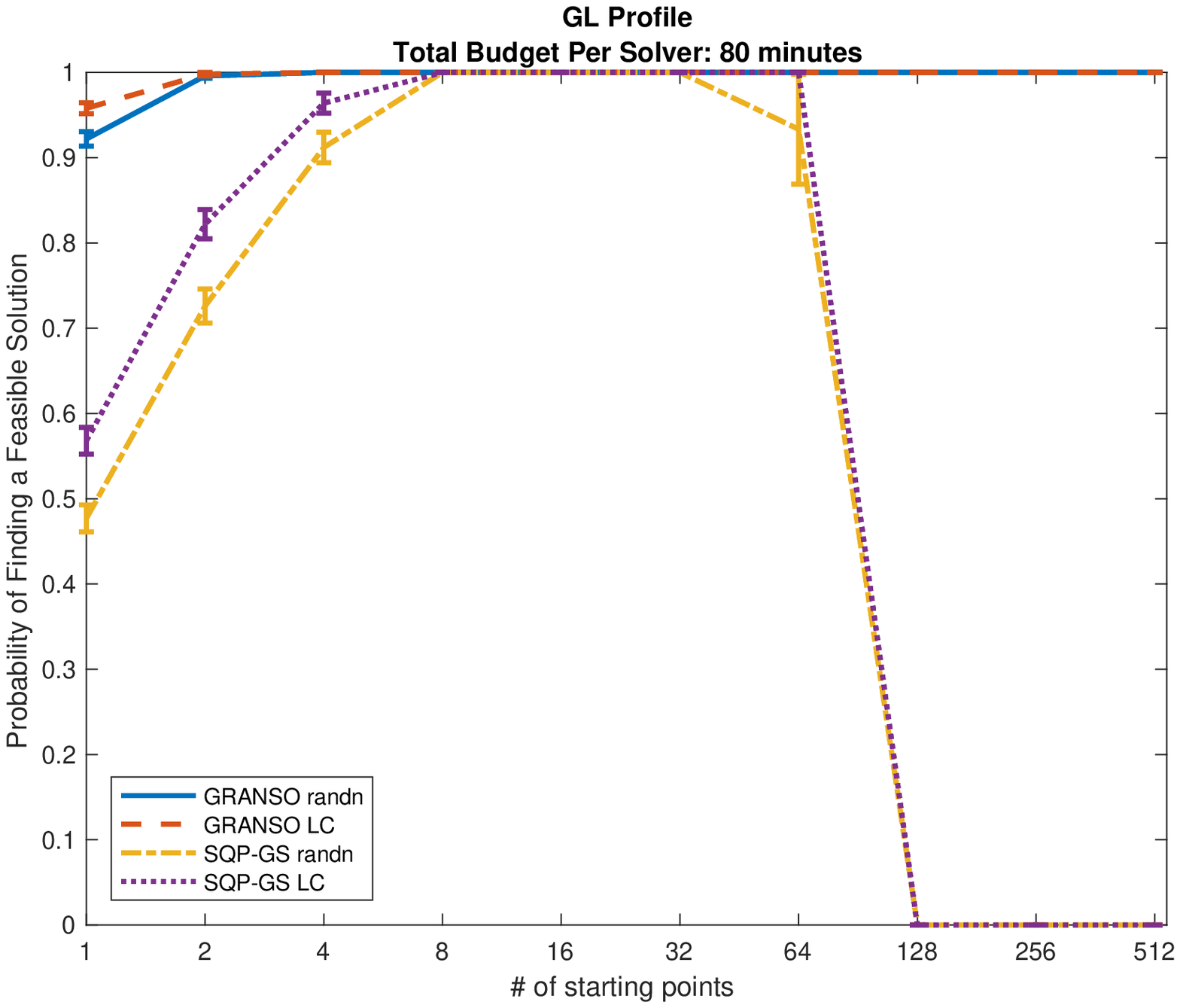}}\qquad
\subfloat{\includegraphics[width=0.40\textwidth]{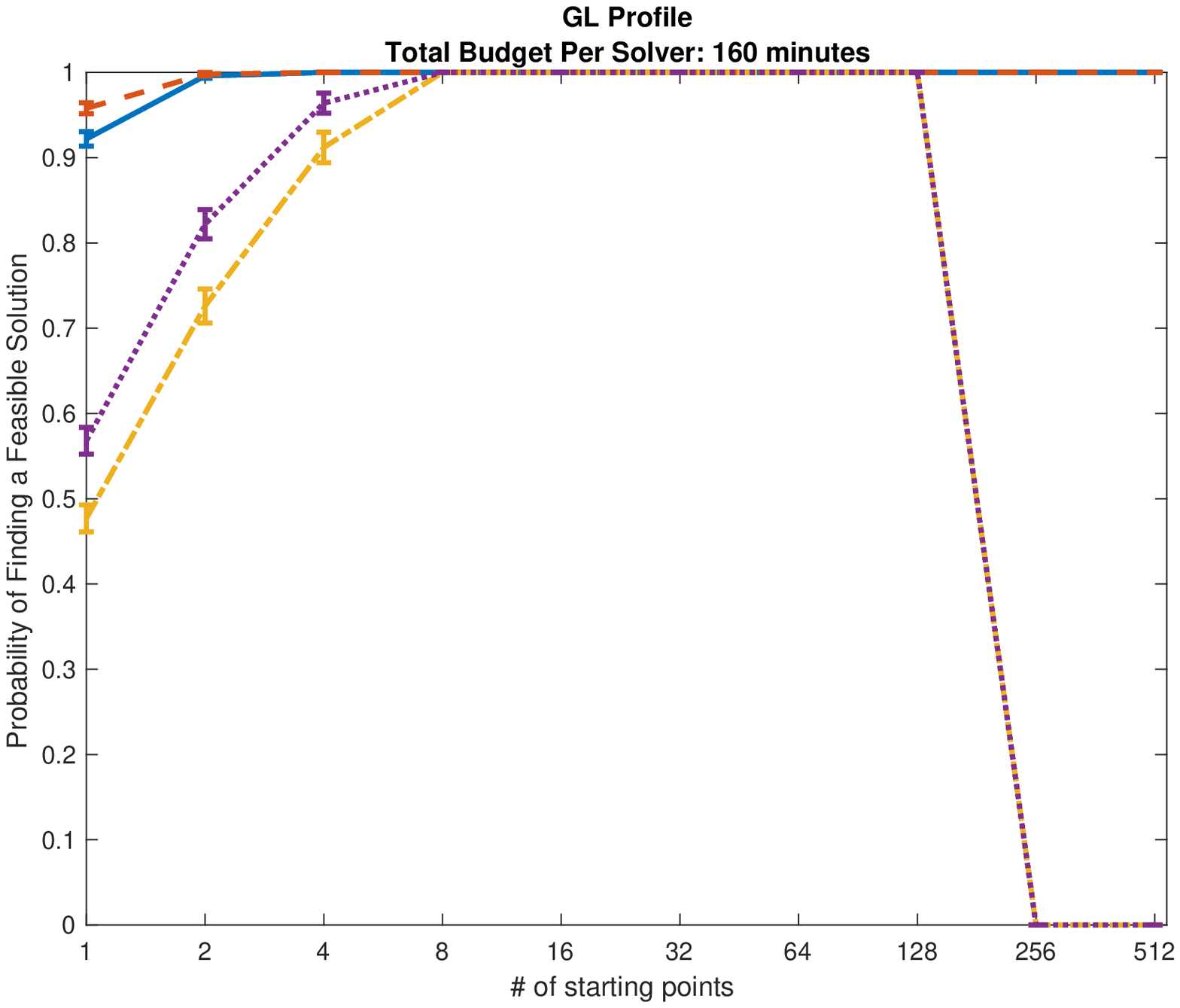}}
\caption{GL-Profiles for probability of finding at least one feasible solution.}\label{fig:gltwo}
\end{figure}

Similar to our $\beta$-RMP benchmark, we plotted GL-Profile curves for each pairing of the two solvers and the two starting sets,
and using the same tolerances to determine whether an iterate is feasible or not.
We created two different panels, one for the expected objective value and one for the expectation of just finding feasible solutions,
respectively shown in Figures~\ref{fig:glone} and \ref{fig:gltwo}.
For each panel, we created four separate GL-Profiles for the following four different budgets: 
20, 40, 80, and 160 minutes of wall-clock time.
Each GL-Profile divides their respective budget over $2^r$ starting points, for $r = 1,\ldots,9$.
Since our initial points were in no particular order, we simply constructed $P$ sets by assigning the first 
$2^r$ points to $S_1$, the next to $S_2$, etc.  Since $N=1000$ is not a power of two,
the less than $2^r$ remaining initial points were simply omitted from each $P$.
Of course, as $r$ increases, $|P|$ decreases and for $r=9$, $P$ will contain only one selection, namely the 
first 512 starting points.

In Figure~\ref{fig:glone}, our GL-Profiles consider the expected best objective value to be obtained on the feasible set.
For a small budget of just 20 minutes total, we see that GRANSO initialized from 32 of the LC points is the best strategy.
SQP-GS is so slow that it actually does not make sense to run more than one single point, since even running from just two points for 10 minutes each returns a worse result.
But with the total budget increased to 40 minutes, we see that running SQP-GS from two or four starting is better than just one.
Still, from the \texttt{randn} test set, GRANSO provides a better optimal strategy for budget allocation than SQP-GS does.
However, on the LC test set, we see that optimal allocations for GRANSO and SQP-GS offer similar solution quality.
For budgets of 80 minutes or longer, we now see that SQP-GS initialized from four to eight LC points provides the best strategy.
This slightly larger sample size significantly increases the probability that one of the starting points might be from 
the right side of the plots shown in Figures~\ref{fig-randbestvals} and \ref{fig-ineqbestvals},
where SQP-GS generally found a larger percentage of higher quality solutions than GRANSO did.

In Figure~\ref{fig:gltwo}, the GL-Profiles only consider the probability of finding the feasible region.
Here, the scenario is quite different compared to the objective value criterion. 
Now we see that GRANSO is always better than SQP-GS at finding the feasible set, regardless of budget.
This is of course not surprising given feasible/infeasible distributions shown in Figures~\ref{fig-randbestvals} and \ref{fig-ineqbestvals}.

\section{Concluding Remarks}
\label{sec:conclusion}
In our evaluation of using nonsmooth, nonconvex, constrained optimization solvers, namely SQP-GS and GRANSO, 
to efficiently design input shapers that effectively damp undesirable oscillatory modes of mechanical systems, 
we have observed noticeable performance differences between them, some expected and some surprising.
At a high-level, and as we would expect, SQP-GS was at least a magnitude of order slower than GRANSO,
highlighting that the large computational burden of SQP-GS may simply make it a too impractical choice for many users.
As illustrated by our $\beta$-RMP benchmark, 
even when SQP-GS was given up to twelve times longer to run than the median running time for GRANSO on the LC starting points,
SQP-GS was still not at all competitive with GRANSO in terms of minimizing the objective on the feasible set.
More surprising is that SQP-GS was not particularly effective at finding the feasible set at all, regardless of the budget, 
having a success rate of only about 50\%.  Meanwhile, GRANSO was particularly adept at finding the feasible set, rarely 
failing to do so.  
Nevertheless, SQP-GS still managed to return the majority of the best feasible solutions, even though GRANSO found the very best solution.
From this limited perspective, one might be hard pressed to choose which of the two solvers would be better for the task of designing input shapers.

However, our newly introduced GL-Profiles do shed further light on how to select an optimal strategy that maximizes the likelihood of getting 
the best quality solution to a given problem under a fixed computational budget.
For those with who either need solutions quickly, or have limited computational resources, 
GRANSO started with 32 or more of the LC starting points provides the best strategy.
Otherwise, for users that are less restricted in their computational budget 
(specifically, at least 160 minutes total compute time on fast modern hardware),
the best strategy is to use SQP-GS, dividing the budget equally amongst from eight or so initial points,
constructed to satisfy as many of the constraints as possible (i.e. the LC test set).
Thus, with aid of both $\beta$-RMPs and our new GL-Profiles, 
budget-dependent best practices for designing input shapers via nonsmooth optimization solvers are apparent.

\section{Acknowledgement}
The presented research was supported by the Czech Science Foundation under the project 16-17398S.

\bibliographystyle{alpha}
\bibliography{references}
\end{document}